\documentclass[12pt]{article}
\usepackage{amssymb}
\usepackage{amsmath}
\usepackage{amsfonts}

\newcommand\qed{\hfill $\square$}

\newcommand{\RR}{\mathbb R}

\newtheorem{thm}{Theorem}
\newtheorem{prop}[thm]{Proposition}
\newtheorem{lem}[thm]{Lemma}
\newtheorem{cor}[thm]{Corollary}
\newtheorem{rem}[thm]{Remark}

\newcommand{\beqn}{\begin{equation}}
\newcommand{\eeqn}{\end{equation}}
\newcommand{\bear}{\begin{eqnarray}}
\newcommand{\eear}{\end{eqnarray}}
\newcommand{\bean}{\begin{eqnarray*}}
\newcommand{\eean}{\end{eqnarray*}}

\begin{document}

\begin{center}
\bf GRADIENT ESTIMATES FOR A DEGENERATE PARABOLIC EQUATION
WITH GRADIENT ABSORPTION AND APPLICATIONS \rm
\end{center}

\vspace{0.3cm}

\centerline{Jean-Philippe Bartier\footnote{CEREMADE, Universit\'e 
Paris-Dauphine, Place du Mar\'echal de Lattre de Tassigny, F--75775
Paris Cedex 16, France. {\tt E-mail: bartier@ceremade.dauphine.fr}} and
Philippe Lauren\c cot\footnote{Institut de
Math\'ematiques de Toulouse, CNRS UMR~5219, Universit\'e Paul Sabatier
(Toulouse~III), 118 route de Narbonne, F--31062 Toulouse Cedex 9,
France. {\tt E-mail: laurenco@mip.ups-tlse.fr}}} 

\vspace{0.2cm}

\begin{abstract}
Qualitative properties of non-negative solutions to a quasilinear degenerate
parabolic equation with an absorption term depending solely on the
gradient are shown, providing information on the competition between
the nonlinear diffusion and the nonlinear absorption. In particular,
the limit as $t\to\infty$ of the $L^1$-norm of integrable solutions is
identified, together with the rate of expansion of the support for
compactly supported initial data. The persistence of dead cores is
also shown. The proof of these results strongly relies on gradient
estimates which are first established. 
\end{abstract}

\vspace{0.2cm}

\section{Introduction}\label{int}

\setcounter{thm}{0}
\setcounter{equation}{0}

We investigate the properties of non-negative and bounded continuous
solutions to the Cauchy problem  
\bear
\label{vhj1}
\partial_t u - \Delta_p u + \vert\nabla u\vert^q & = & 0\ \,, 
\quad (t,x)\in Q_\infty:=(0,\infty)\times\RR^N\,,\\
\label{vhj2}
u(0) & = & u_0\ge 0\,, \quad x\in\RR^N\,,
\eear
the parameters $p$ and $q$ ranging in $(2,\infty)$ and $(1,\infty)$,
respectively, and the $p$-Laplacian operator $\Delta_p$ being defined by 
$$
\Delta_p u := \mbox{ div } \left( \vert\nabla u\vert^{p-2}\ \nabla u
\right)\,.
$$
When $p>2$, (\ref{vhj1}) is a quasilinear degenerate parabolic
equation with a nonlinear absorption term $|\nabla u|^q$ depending
solely on the gradient of $u$, and reduces to the semilinear diffusive
Hamilton-Jacobi equation
\beqn
\label{vhj3}
\partial_t v - \Delta v + \vert\nabla v\vert^q = 0\ \;\;\mbox{ in
}\;\; Q_\infty\,,
\eeqn
when $p=2$. Several recent papers have been devoted to the study of
properties of non-negative solutions to (\ref{vhj3}) with a particular
emphasis on the large time behaviour which turns out to depend
strongly on the value of the parameter $q\in (0,\infty)$
\cite{ATU04,BKL04,BL99,BLSS02,BRV96,BRV97,Gi05}. 

One of the keystones of these investigations are optimal gradient
estimates of the form $\|\nabla\left( v^\alpha \right)(t)\|_\infty\le
C(\|v(0)\|_\infty)\ t^{-\beta}$ for suitable exponents $\alpha\in
(0,1)$ and $\beta>0$, both depending on $N$ and $q$
\cite{BL99,GGK03}. Not only do such estimates provide an instantaneous
smoothing effect from $L^\infty(\RR^N)$ to $W^{1,\infty}(\RR^N)$ but
temporal decay estimates as well, the latter being the starting point
of a precise study of the large time dynamics. Let us recall here that
the proof of the above-mentioned gradient estimates relies on a 
modification of the Berstein technique \cite{BL99,GGK03}.

\medskip

Owing to the nonlinearity of the diffusion term when $p>2$, the
availability of similar gradient estimates for solutions to (\ref{vhj1}), 
(\ref{vhj2}) is unclear and is actually
our first result. More precisely, for $p>2$ and $q>1$, we introduce
the exponents $\alpha_p\in (0,1)$ and $\beta_{p,q}\in (0,1)$ defined by 
\beqn
\label{spirou}
\frac{1}{\alpha_p}  := \frac{p-1}{p-2} - \frac{N-1}{p(N+3)-2(N+1)}
\;\;\mbox{ and }\;\; \beta_{p,q} := \max{\left\{ \alpha_p ,
\frac{q-1}{q} \right\}}\,. 
\eeqn

\begin{thm}\label{tha1}
Consider a non-negative initial condition $u_0\in\mathcal{BC}(\RR^N)$.
There is a non-negative viscosity solution
$u\in\mathcal{BC}([0,\infty)\times \RR^N)$ to (\ref{vhj1}),
(\ref{vhj2}) such that 
\beqn
\label{se1}
0\le u(t,x)\le \|u_0\|_\infty\,, \quad (t,x)\in Q_\infty\,,
\eeqn 
\bear
\label{ge1}
\left\vert \nabla\left( u^{\alpha_p} \right)(t,x) \right\vert & \le & 
C(p,N)\ \Vert u(s)\Vert_\infty^{(p \alpha_p+2-p)/p}\ (t-s)^{-1/p}\,,
\\
\nonumber
 & & \\
\label{ge2}
\left\vert \nabla\left( u^{\beta_{p,q}} \right)(t,x) \right\vert & \le
& C(p,q,N)\ \Vert u(s)\Vert_\infty^{(q \beta_{p,q} +1- q)/q}\
(t-s)^{-1/q}\,, 
\eear
and
\beqn
\label{weakf}
\int_{\RR^N} (u(t,x)-u(s,x))\ \vartheta(x)\ dx + \int_s^t \int_{\RR^N}
\left( |\nabla u|^{p-2} \nabla u \cdot \nabla\vartheta + |\nabla u|^q\
\vartheta \right)\ dxd\tau = 0
\eeqn
for $t>s\ge 0$ and $\vartheta\in\mathcal{C}_0^\infty(\RR^N)$.

Furthermore, this solution is unique if $u_0\in\mathcal{BUC}(\RR^N)$.
\end{thm}

Let us emphasize that the main contribution of Theorem~\ref{tha1} is
the estimates (\ref{ge1}), (\ref{ge2}), and not the existence of a
viscosity solution to (\ref{vhj1}) which could probably be obtained by
alternative approaches. But, owing to the poor regularity of the
solutions to (\ref{vhj1}), (\ref{vhj2}), we cannot prove (\ref{ge1})
and (\ref{ge2}) directly and instead use an approximation
procedure. Indeed, the proof of (\ref{ge1}) and (\ref{ge2}) relies on
a modification of the Bernstein technique. It requires the study of
the partial differential equation solved by $|\nabla\varphi(u)|^2$ for
a suitably chosen function $\varphi$ and thus some regularity which is
not available for solutions to (\ref{vhj1}), (\ref{vhj2}). The existence
part of Theorem~\ref{tha1} is in fact an intermediate step in the
proof of (\ref{ge1}) and (\ref{ge2}).  

It is clear from (\ref{ge1}) and (\ref{ge2}) with $s=0$ that they lead
to different temporal decay estimates. In fact, as we shall see below,
(\ref{ge1}) results from the diffusive part of (\ref{vhj1}) while
(\ref{ge2}) stems from the absorption term. In particular, it is
worth mentioning that (\ref{ge1}) is also valid for
non-negative solutions to the $p$-Laplacian equation 
\beqn
\label{ple}
\partial_t w - \Delta_p w = 0 \;\;\mbox{ in }\;\; Q_\infty\,,
\eeqn
 which seems to be new for $N\ge
2$. When $N=1$, it has been proved in \cite[Theorem~2]{EV86}. Also,
(\ref{ge2}) is true for non-negative viscosity solutions to the
Hamilton-Jacobi equation 
\beqn
\label{hje}
\partial_t h +|\nabla h|^q=0 \;\;\mbox{ in }\;\; Q_\infty\,,
\eeqn
and can be deduced from \cite[Theorem~I.1]{Li85}. For $p=2$, similar
gradient estimates have been obtained in \cite{BL99,GGK03} with
$\alpha_2=\beta_{2,q}=(q-1)/q$.  

The previous gradient estimates may be improved for non-negative,
radially symmetric, and non-increasing initial data. 

\begin{thm}\label{tha2}
Assume that the initial condition $u_0\in\mathcal{BC}(\RR^N)$ is
non-negative, radially symmetric, and non-increasing. There is a
non-negative viscosity solution $u$ to (\ref{vhj1}), (\ref{vhj2})
satisfying (\ref{se1}), (\ref{weakf}) and such that 
$$
x\longmapsto u(t,x) \;\;\mbox{ is non-negative, radially symmetric,
and non-increasing, }
$$  
\bear
\label{ge1r}
\left\vert \nabla\left( u^{(p-2)/(p-1)} \right)(t,x) \right\vert & \le
& C(p,N)\ \Vert u(s)\Vert_\infty^{(p-2)/p(p-1)}\ (t-s)^{-1/p}\,,\\ 
\nonumber
& & \\
\label{ge2ra}
\left\vert \nabla\left( u^{(q-1)/q} \right)(t,x) \right\vert & \le & 
\frac{(q-1)^{(q-1)/q}}{q}\ t^{-1/q} \;\;\mbox{ if }\;\; q\ge p-1\,, 
\eear
and
\beqn
\label{ge2rb}
\left\vert \nabla\left( u^{(p-2)/(p-1)} \right)(t,x) \right\vert \le
C(p,q)\ \Vert u(s)\Vert_\infty^{(p-1-q)/q(p-1)}\ (t-s)^{-1/q}
\;\;\mbox{ if }\;\; q\in (1,p-1)\,, 
\eeqn
for $t>s\ge 0$.
\end{thm} 

Theorem~\ref{tha2} is proved as Theorem~\ref{tha1} for $N=1$. We 
will thus only give the proof of the latter.

Here again, the gradient estimate (\ref{ge1r}) is valid for
non-negative solutions to the $p$-Laplacian equation (\ref{ple}) with radially
symmetric and non-increasing initial data and is easily seen to be
optimal in that case: indeed, the Barenblatt solution to the
$p$-Laplacian equation (\ref{ple}) is given by  
$$
\mathcal{B}(t,x) = t^{-N\eta}\ \left( 1 - \gamma_p \left( \frac{\vert
x\vert}{t^\eta} \right)^{p/(p-1)} \right)_+^{(p-1)/(p-2)}\,, \quad
(t,x)\in (0,\infty)\times\RR^N\,, 
$$
(see, e.g., \cite[Ch.~XI, Eq.~(1.6)]{DB93}) and $\nabla \left( 
\mathcal{B}^\vartheta \right)(t,x)$ is bounded only
for $\vartheta\ge (p-2)/(p-1)$. 

\begin{rem}\label{recomp}
Since we are mainly interested in qualitative properties of solutions
to (\ref{vhj1}), (\ref{vhj2}), we leave aside the question of
uniqueness of such solutions for initial data in
$\mathcal{BC}(\RR^N)\setminus\mathcal{BUC}(\RR^N)$. Nevertheless,
since the solutions in Theorems~\ref{tha1} and~\ref{tha2}
are constructed as limits of classical solutions, they still enjoy a
comparison principle. More precisely, if $u_0$ and $\hat{u}_0$ are two
non-negative functions in $\mathcal{BC}(\RR^N)$ such that
$u_0\le\hat{u}_0$, then the corresponding solutions $u$ and $\hat{u}$
to (\ref{vhj1}) with initial data $u_0$ and $\hat{u}_0$ constructed in
Theorem~\ref{tha1} satisfy $u(t,x)\le \hat{u}(t,x)$ for all $(t,x)\in
Q_\infty$. This fact will be used repeatedly in the sequel. 
\end{rem} 

Several qualitative properties follow from the previous gradient
estimates. As a first consequence, we derive temporal decay estimates
in $W^{1,\infty}(\RR^N)$ for non-negative and integrable
solutions to (\ref{vhj1}), (\ref{vhj2}). We set 
\beqn
\label{spip}
q_*:=p-\frac{N}{N+1}\,, \quad \xi:= \frac{1}{q(N+1)-N}\,, \quad \eta
:= \frac{1}{N(p-2)+p}\,. 
\eeqn

\begin{prop}\label{pra3}
Assume that 
\beqn
\label{fantasio}
u_0\in L^1(\RR^N)\cap\mathcal{BC}(\RR^N)\,, \;\; u_0\ge 0\,,
\eeqn
and denote by $u$ the corresponding viscosity solution to
(\ref{vhj1}), (\ref{vhj2}) constructed in Theorem~\ref{tha1}. Then
$u\in \mathcal{C}([0,\infty);L^1(\RR^N))$. 

Let $t>0$. If $q\in
(1,q_*)$, then 
\bear
\label{dec1}
\Vert u(t)\Vert_\infty & \le & C\ \Vert u_0\Vert_1^{q\xi}\
t^{-N\xi}\,, \\ 
\label{dec2}
\Vert\nabla u(t)\Vert_\infty & \le & C\ \Vert u_0\Vert_1^{\xi}\
t^{-(N+1)\xi}\,, 
\eear
while, if $q>q_*$, 
\bear
\label{dec3}
\Vert u(t)\Vert_\infty & \le & C\ \Vert u_0\Vert_1^{p\eta}\
t^{-N\eta}\,, \\ 
\label{dec4}
\Vert\nabla u(t)\Vert_\infty & \le & C\ \Vert u_0\Vert_1^{2\eta}\
t^{-(N+1)\eta}\,. 
\eear
\end{prop}

Recall that the $L^\infty$-norm of non-negative and integrable
solutions $w$ to the $p$-Laplacian equation (\ref{ple}) decays as
$t^{-N\eta}$ 
\cite[Theorem~3]{HV81}. However this decay might be enhanced by the
nonlinear absorption term and this is indeed the case for $q\in
(1,q_*)$. Indeed, $t^{-N\xi}\le t^{-N\eta}$ for $t\ge 1$ and $q\in
(1,q_*)$. According to Proposition~\ref{pra3}, we thus expect the
nonlinear absorption term to be negligible as $t\to\infty$ for $q>q_*$
and the large time dynamics to feel the effects of the absorption only
for $q\in (1,q_*)$. The next result is a further step in that
direction.

\medskip

It readily follows from (\ref{vhj1}) and the non-negativity of $u$
that $t\longmapsto \Vert u(t)\Vert_1$ is a non-increasing and
non-negative function. Introducing 
\beqn
\label{gaston}
I_1(\infty) := \lim_{t\to\infty} \Vert u(t)\Vert_1 = \inf_{t\ge
0}{\left\{ \Vert u(t)\Vert_1 \right\}} \in \left[0,\Vert u_0\Vert_1
\right] \,, 
\eeqn
we study the possible values of $I_1(\infty)$. 

\begin{prop}
\label{pra4}
Assume that $u_0$ satisfies (\ref{fantasio}) with $\Vert u_0\Vert_1>0$
and denote by $u$ the
corresponding viscosity solution to (\ref{vhj1}), (\ref{vhj2})
constructed in Theorem~\ref{tha1}. Then $I_1(\infty)>0$ if and
only if $q>q_*$, the parameter $q_*$ being defined in (\ref{spip}).  
\end{prop}

Since  $\|w(t)\|_1=\|w(0)\|_1$ for all $t\ge 0$ for non-negative and 
integrable solutions $w$ to the $p$-Laplacian equation (\ref{ple}), 
we realize that the absorption term is not strong
enough for $q>q_*$ to drive the $L^1$-norm of $u(t)$ to zero as
$t\to\infty$, thus indicating a diffusion-dominated behaviour for 
large times. For $q\in (p-1,p)$ Proposition~\ref{pra4} is already 
proved in \cite[Theorems~1.3 \&~1.4]{ATU04} by a different method.  

\medskip

We next turn to a property which marks a striking difference between
the semilinear case $p=2$ and the quasilinear case $p>2$ corresponding
to \textit{slow diffusion}, namely the finite speed of
propagation. Since the support of non-negative and compactly supported
solutions $w$ to the $p$-Laplacian equation (\ref{ple}) grows as
$t^\eta$, it is natural to wonder whether the absorption term will
slow down this process.  

\begin{thm}\label{tha5}
Assume that $u_0$ fulfils (\ref{fantasio}) and is compactly
supported, and denote by $u$ the corresponding solution to
(\ref{vhj1}), (\ref{vhj2}). For $t\ge 0$ we put 
\beqn
\label{rayon}
\varrho(t) := \inf{\left\{ R>0 \;\;\mbox{ such that }\;\; u(t,x)=0
\;\;\mbox{ for }\;\; \vert x\vert>R \right\}}\,. 
\eeqn
Then $\varrho(t)<\infty$ for all $t\ge 0$ and: 
\begin{itemize}
\item[(i)] If $q\in (1,p-1)$ then
\beqn
\label{rad1}
\limsup_{t\to\infty} \varrho(t) < \infty\,.
\eeqn
\item[(ii)] If $q=p-1$ then
\beqn
\label{rad2}
\varrho(t) \le C\ (1 + \ln{t} ) \;\;\mbox{ for }\;\; t\ge 1\,.
\eeqn
\item[(iii)] If $q\in (p-1,q_*)$ then
\beqn
\label{rad3}
\varrho(t) \le C\ t^{(q-p+1)/(2q-p)} \;\;\mbox{ for }\;\; t\ge 1\,.
\eeqn
\item[(iv)] If $q\ge q_*$ then
\beqn
\label{rad4}
\varrho(t) \le C\ t^\eta \;\;\mbox{ for }\;\; t\ge 1\,.
\eeqn
\end{itemize}
\end{thm}

Here again, the absorption term seems to have no real effect on the
expansion on the support of $u(t)$ for $q>q_*$ as the upper bound
(\ref{rad4}) is exactly the growth rate 
of the support for non-negative and
compactly supported solutions $w$ to the $p$-Laplacian equation
(\ref{ple}). But, as soon as $q$ is below $q_*$, the dynamics starts
to feel the effects of the absorption term and the expansion
of the support of $u(t)$ slows down. It even stops for $q\in (1,p-1)$.
In that case, the support of $u(t)$ remains \textit{localized} in a 
fixed ball of $\RR^N$: such a property is already enjoyed by compactly
supported non-negative solutions to second-order degenerate parabolic
equations with a sufficiently strong absorption involving the solution
only as, for instance, $\partial_t z - \Delta_p z
+ z^r = 0$ in $Q_\infty$ when $r\in (1,p-1)$ \cite{DV85,Ka87,Yu96}. It has
apparently remained unnoticed for second-order degenerate parabolic
equations with an absorption term depending solely on the gradient. In
our case, this property is clearly reminiscent of that enjoyed by the
solutions $h$ to the Hamilton-Jacobi equation (\ref{hje}): namely, the
support of $h(t)$ does not evolve through time evolution
\cite{Bl85}. Finally, for $q\in (p-1,q_*)$, compactly supported
self-similar solutions to (\ref{vhj1}) are constructed and the
boundaries of their support evolve at the speed given by the
right-hand side of (\ref{rad3}). 

\medskip

As a by-product of the proof of Theorem~\ref{tha5} we obtain improved
decay estimates for the $L^1$-norm of solutions to (\ref{vhj1}),
(\ref{vhj2}) with compactly supported initial data. 

\begin{cor}\label{cora6}
Assume that $u_0$ fulfils (\ref{fantasio}) and is compactly
supported. Then 
\begin{itemize}
\item[(i)] If $q\in (1,p-1)$ then
\beqn
\label{decL11}
\Vert u(t)\Vert_1 \le C\ t^{-1/(q-1)}\,, \quad t\ge 2\,.
\eeqn
\item[(ii)] If $q=p-1$ then
\beqn
\label{decL12}
\Vert u(t)\Vert_1 \le C\ t^{-1/(q-1)}\ \left( \ln{t}
\right)^{1/\xi(q-1)} \;\;\mbox{ for }\;\; t\ge 2\,. 
\eeqn
\item[(iii)] If $q\in (p-1,q_*)$ then
\beqn
\label{decL13}
\Vert u(t)\Vert_1 \le C\ t^{-((N+1)(q_*-q))/(2q-p)} \;\;\mbox{ for
}\;\; t\ge 2\,. 
\eeqn
\item[(iv)] If $q=q_*$ then
\beqn
\label{decL14}
\Vert u(t)\Vert_1 \le C\ \left( \ln{t} \right)^{-1/(q-1)} \;\;\mbox{
for }\;\; t\ge 2\,. 
\eeqn
\end{itemize}
\end{cor}

For $q\in (p-1,q_*]$, Theorem~\ref{tha5} and Corollary~\ref{cora6} are
already proved in \cite[Theorems~1.1 \&~1.2]{ATU04} by a completely
different approach. In addition, for non-compactly supported initial data, 
temporal decay estimates involving the behaviour of $u_0$ for large
values of $x$ are obtained in \cite[Theorem~1.3]{ATU04} for the
$L^1$-norm of $u$. Let us also mention that the decay rate of $\Vert
u(t)\Vert_1$ for $q\in (1,p-1)$ is the same as the one obtained in
\cite{Bl85} for non-negative and compactly supported solutions to the
Hamilton-Jacobi equation (\ref{hje}). The bound (\ref{decL11}) then
provides another clue of the dominance of the absorption term for
$q\in (1,p-1)$. That it is indeed true is shown in \cite{LVxx}. 

\medskip 

For $q\in (1,p-1)$, it follows from Theorem~\ref{tha5}~(i) that the support 
of the solutions to 
(\ref{vhj1}), (\ref{vhj2}) with compactly supported initial data remains 
bounded through time evolution. A natural counterpart of this phenomenon 
is to study what happens to a solution to (\ref{vhj1}), (\ref{vhj2}) starting 
from an initial condition vanishing inside a ball of $\RR^N$. It turns out that, 
if the radius of the ball is sufficiently large, the solution still vanishes 
inside of a smaller ball for all times, a phenomenon which may be called the 
persistence of \textit{dead cores}. 

\begin{prop}\label{pra7}
Consider a non-negative initial condition $u_0\in\mathcal{BC}(\RR^N)$ such that
\beqn
\label{pdc1}
u_0(x) = 0 \;\;\mbox{ if }\;\; \vert x\vert\le R_0
\eeqn
for some $R_0>0$, and denote by $u$ the corresponding solution to
(\ref{vhj1}), (\ref{vhj2}) constructed in Theorem~\ref{tha1}. If $q\in
(1,p-1)$ there is a constant $\delta_0=\delta_0(p,q)>0$ such that, if
$R_0\ge \delta_0\ \Vert u_0\Vert_\infty^{(p-1-q)/(p-q)}$ then 
$$
u(t,x) = 0 \;\;\mbox{ if }\;\; \vert x\vert \le R_0 - \delta_0\ \Vert
u_0\Vert_\infty^{(p-1-q)/(p-q)} \;\;\mbox{ and }\;\; t\ge 0\,. 
$$
\end{prop}

The proof of Proposition~\ref{pra7} is in fact  quite similar to that of
Theorem~\ref{tha5}~(i).

\medskip

This paper is organized as follows: gradient estimates for an
approximation of (\ref{vhj1}) are established in Section~\ref{ge} by a
modified Bernstein technique with the help of a trick introduced in
\cite{Be81} to obtain gradient estimates for the porous medium
equation. Theorems~\ref{tha1} and~\ref{tha2} are then proved in
Section~\ref{exist}. Sections~\ref{tde} and~\ref{lvl1} are devoted to
integrable initial data for which we prove Propositions~\ref{pra3}
and~\ref{pra4}. We focus on compactly supported initial data in
Section~\ref{csid} where Theorem~\ref{tha5} and Corollary~\ref{cora6}
are proved. The persistence of dead cores is studied in
Section~\ref{pode} while the proof of a technical lemma from
Section~\ref{ge} is postponed to the appendix.  
 

\section{Gradient estimates}\label{ge}

\setcounter{thm}{0}
\setcounter{equation}{0}

As already mentioned the proof of the gradient estimates (\ref{ge1}) 
and (\ref{ge2}) rely on a modified Bernstein technique: owing to the 
degeneracy of the diffusion we cannot expect (\ref{vhj1}) to have 
smooth solutions and we thus need to use an approximation procedure. 
We first report the following technical lemma.

\begin{lem}\label{leb1}
Let $a$ and $b$ be two non-negative functions in
$\mathcal{C}^2([0,\infty))$ and $u$ be a classical solution to 
\beqn
\label{b1}
\partial_t u - \mbox{ div } \left( a\left( \vert\nabla u\vert^2
\right)\ \nabla u \right) + b\left( \vert\nabla u\vert^2 \right) = 0\
\;\;\mbox{ in }\;\; Q_\infty\,.
\eeqn
Consider next a $\mathcal{C}^3$-smooth increasing function $\varphi$ 
and set $v:= \varphi^{-1}(u)$ and $w:=\vert\nabla v\vert^2$. Then $w$
satisfies the following differential inequality 
\beqn
\label{b2}
\partial_t w - \mathcal{A} w - \mathcal{V}\cdot
\nabla w + 2\ \mathcal{R}_1\ w^2 + 2\ \mathcal{R}_2\ w \le 0 \;\;\mbox{ in
}\;\; Q_\infty\,,
\eeqn
where $\mathcal{A}$, $\mathcal{R}_1$ and $\mathcal{R}_2$ are given by 
\beqn
\label{b3}
\mathcal{A}w := a\ \Delta w + 2 a'\ (\nabla u)^t D^2 w \nabla u\,,
\eeqn
\beqn
\label{b4}
\mathcal{R}_1 := - a\ \left( \frac{\varphi''}{\varphi'}
\right)' - \left( (N-1)\ \frac{{a'}^2}{a} + 4\ a'' \right)\ \left(
\varphi' \varphi'' \right)^2\ w^2 - 2\ a'\ w \left( 2 {\varphi''}^2 +
\varphi' \varphi''' \right)\,,
\eeqn
\beqn
\label{b5}
\mathcal{R}_2 := \frac{\varphi''}{{\varphi'}^2}\ \left( 2\ b'\
{\varphi'}^2\ w - b \right)\,,
\eeqn
while $\mathcal{V}$ is given by (\ref{ap1}) below. Here and in the
following we omit the variable in $a$, $b$ and $\varphi$ and their
derivatives. 

Furthermore, if $\varphi$ is convex, $a$ is non-decreasing and
$x\longmapsto u(t,x)$ is radially symmetric and non-increasing for
each $t\ge 0$, then $\mathcal{R}_1$ may be replaced by
$\mathcal{R}_1^r$ given by 
\beqn
\label{b6}
\mathcal{R}_1^r := - a\ \left( \frac{\varphi''}{\varphi'} \right)' -
4\ a''\ \left( \varphi' \varphi'' \right)^2\ w^2 - 2\ a'\ w \left( 2
{\varphi''}^2 + \varphi' \varphi''' \right)\,,
\eeqn
\end{lem}

The proof of Lemma~\ref{leb1} is rather technical and is postponed to
the appendix. We however emphasize that it uses a trick introduced by
B\'enilan \cite{Be81} to prove gradient estimates for solutions to the
porous medium equation in several space dimensions. It is also worth 
noticing that $\mathcal{R}_1=\mathcal{R}_1^r$ for $N=1$. 

\medskip

Consider next a non-negative function $u_0\in \mathcal{BC}(\RR^N)$. There
is a sequence of functions $(u_{0,k})_{k\ge 1}$ such that, for each
integer $k\ge 1$, $u_{0,k}\in\mathcal{BC}^\infty(\RR^N)$,
\beqn
\label{ex1}
0 \le u_{0,k}(x) \le u_{0,k+1}(x) \le u_0(x)\,, \quad x\in\RR^N\,,
\eeqn
and $(u_{0,k})$ converges uniformly towards $u_0$ on compact subsets
of $\RR^N$. In addition, if $u_0\in W^{1,\infty}(\RR^N)$ we may assume
that
\beqn
\label{ex2}
\Vert\nabla u_{0,k}\Vert_\infty \le \left( 1 + \frac{K_1}{k} \right)\
\Vert \nabla u_0\Vert_\infty\,,
\eeqn
for some constant $K_1>0$ depending only on the approximation process.
Next, since $\xi\longmapsto \vert\xi\vert^{p-2}$ and $\xi\longmapsto
\vert\xi\vert^q$ are not regular enough for small
values of $p$ and $q$, we set
\beqn
\label{ex3}
a_\varepsilon(\xi) := \left( \varepsilon^2 + \xi \right)^{(p-2)/2}
\;\;\mbox{ and }\;\; b_\varepsilon(\xi) := \left( \varepsilon^2 + \xi
\right)^{q/2} - \varepsilon^q\,, \quad \xi\ge 0\,,
\eeqn
for $\varepsilon\in (0,1/2)$. Then, given 
\beqn
\label{ex3b}
0 < \gamma \le \min{\left\{ \frac{3}{4} , 2\beta_{p,q} , q ,
\frac{q+2}{2} \right\}} \,,
\eeqn
the Cauchy problem
\bear
\label{ex4}
\partial_t u_{k,\varepsilon} - \mbox{ div } \left( a_\varepsilon\left(
\vert\nabla u_{k,\varepsilon} \vert^2 \right)\ \nabla
u_{k,\varepsilon} \right) + b_\varepsilon\left( \vert\nabla
u_{k,\varepsilon}\vert^2 \right) & = & 0\,,
\quad (t,x)\in Q_\infty\,,\\
\label{ex5} 
u_{k,\varepsilon}(0) & = & u_{0,k}+\varepsilon^\gamma \,, \quad
x\in\RR^N\,, 
\eear 
has a unique classical solution $u_{k,\varepsilon} \in
\mathcal{C}^{(3+\delta)/2,3+\delta}([0,\infty)\times\RR^N)$ for some
$\delta\in (0,1)$ \cite{LSU88}. Observing that $\varepsilon^\gamma$ and
$\Vert u_0\Vert_\infty + \varepsilon^\gamma$ are solutions to (\ref{ex4})
with $\varepsilon^\gamma \le u_{k,\varepsilon}(0,x) \le \Vert
u_0\Vert_\infty + \varepsilon^\gamma$, the comparison principle warrants
that
\beqn
\label{ex4b}
\varepsilon^\gamma \le u_{k,\varepsilon}(t,x) \le \Vert u_0\Vert_\infty +
\varepsilon^\gamma\,, \quad (t,x)\in [0,\infty)\times\RR^N\,.
\eeqn

\medskip

We now turn to estimates on the gradient of $u_{k,\varepsilon}$ and
first point out that, thanks to the regularity of $a_\varepsilon$,
$b_\varepsilon$ and $u_{k,\varepsilon}$, we may use Lemma~\ref{leb1}.
We first take $\varphi(r)=\varphi_0(r):=r$ for $r\ge 0$ so that 
$w=\vert\nabla u_{k,\varepsilon}\vert^2$ and
$\mathcal{R}_1=\mathcal{R}_2=0$. Therefore $w$ satisfies
$$
\partial_t w - \mathcal{A} w - \mathcal{V}\cdot \nabla w \le 0
\;\;\mbox{ in }\;\; Q_\infty\,.
$$
Since $w(0)\le \left\Vert \nabla u_{0,k} \right\Vert_\infty^2$ the
comparison principle ensures that 
\beqn
\label{ex6} 
\left\Vert \nabla u_{k,\varepsilon}(t) \right\Vert_\infty \le \left\Vert 
\nabla u_{0,k} \right\Vert_\infty\,, \quad t\ge 0\,. 
\eeqn

We now establish gradient estimates similar to (\ref{ge1})
and (\ref{ge2}) for $u_{k,\varepsilon}$. We first use the specific
choice of $a_\varepsilon$ and $b_\varepsilon$ to compute
$\mathcal{R}_1$ and $\mathcal{R}_2$.

\begin{lem}\label{leb2}
Introducing $g:= \left( \vert\nabla u_{k,\varepsilon}\vert^2 +
\varepsilon^2 \right)^{1/2}$, we have
\beqn
\label{b10}
\mathcal{R}_1 = - (p-1)\ g^{p-2}\ \left\{ \left(
\frac{\varphi''}{\varphi'} \right)' + \frac{\alpha_p}{1-\alpha_p}\ \left(
\frac{\varphi''}{\varphi'} \right)^2 \right\} + \varepsilon^2\
\mathcal{R}_{11}
\eeqn
with
\bean
\mathcal{R}_{11} & = & (p-2)\ \left( \frac{\varphi''}{\varphi'}
\right)'\ g^{p-4} + \frac{(p-2)(p(N+3)-2(N+1))}{4}\ \left(
\frac{\varphi''}{\varphi'} \right)^2\ g^{p-4} \\
& + & \frac{(p-2)(p(N+3)-2(N+7))}{4}\ \left( \frac{\varphi''}{\varphi'}
\right)^2\ \left( g^2 - \varepsilon^2 \right)\ g^{p-6}\,,
\eean
and
\beqn
\label{b11}
\mathcal{R}_2 = \frac{\varphi''}{{\varphi'}^2}\ \left\{ (q-1)\ g^q +
\varepsilon^q - q\ \varepsilon^2\ g^{q-2} \right\}\,.
\eeqn
\end{lem}

After these preliminary computations we are in a position to state and
prove the main result of this section. 

\begin{prop}\label{prb3} There are positive real numbers $C=C(p,N)$
and $D_1(k)=D_1(k,p,N)$ such that, for $\varepsilon\in (0,1/2)$,
$x\in\RR^N$, and $t\in \left( 0, \varepsilon^{-1/4} \right)$,   
\beqn
\label{b12}
\left\vert \nabla\left( u_{k,\varepsilon}^{\alpha_p} \right)(t,x)
\right\vert \le C\ \left( 1 + D_1(k)\ \varepsilon^{1/4} \right)^{2/p}\
\left( \Vert u_{0,k}\Vert_\infty + \varepsilon^\gamma \right)^{(p
\alpha_p+2-p)/p}\ t^{-1/p}\,.  
\eeqn

There are a positive real number $D_2(k)=D_2(k,p,q,N)$ and a positive function
$\omega\in\mathcal{C}([0,\infty))$ such that $\omega(\varepsilon)\to
0$ as $\varepsilon\to 0$ and 
\bear
\label{b13}
\left\vert \nabla\left( u_{k,\varepsilon}^{\beta_{p,q}} \right)(t,x) 
\right\vert & \le & \frac{\beta_{p,q}}{(q-1)^{1/q} 
(1-\beta_{p,q})^{1/q}}\ \left( \frac{1}{q} 
+ D_2(k)\ \omega(\varepsilon)^{1/2} \right)^{1/q} \\
\nonumber
& & \times \left( \Vert u_{0,k}\Vert_\infty + \varepsilon^\gamma
\right)^{(q \beta_{p,q}+1-q)/q}\ t^{-1/q}
\eear
for $t\in (\left( 0, \omega(\varepsilon)^{-1/2} \right)$, $x\in\RR^N$,
and $\varepsilon\in (0,\min{\{q-1,1/2\}})$.
\end{prop}

The proof of Proposition~\ref{prb3} relies on suitable choices of the
function $\varphi$ in $\mathcal{R}_1$ and $\mathcal{R}_2$. To motivate
the forthcoming choices, we first note that, if
$\varphi(r)=r^{1/\alpha_p}$, then $\mathcal{R}_1=\varepsilon^2\
\mathcal{R}_{11}$ and (\ref{b12}) will in fact be obtained by choosing
a ``small perturbation'' of $r\mapsto r^{1/\alpha_p}$, namely
$\varphi(r)=\varphi_1(r):= \left( 2Kr-r^2 \right)^{1/\alpha_p}$ for
$K$ sufficiently large. Such a choice has already been employed for
the $p$-Laplacian equation in one space dimension $N=1$ for the same
purpose \cite{EV86}. Next, previous investigations for the case $p=2$
suggest that $\varphi(r)=r^{q/(q-1)}$ is a suitable choice in
$\mathcal{R}_2$ \cite{BL99}. However, with this choice of $\varphi$,
$\mathcal{R}_1$ might give a non-positive contribution according to
the value of $p$ and a suitable choice turns out to be
$\varphi(r)=\varphi_2(r):= \beta_{p,q}\ r^{1/\beta_{p,q}}$.  

\medskip

\noindent\textbf{Proof of Proposition~\ref{prb3}.} We first establish
(\ref{b12}). Consider $\mu>0$ to be specified later and put  
$$
K := \sqrt{1+\mu}\ M^{\alpha_p}\,, \quad M:= \Vert u_{0,k}\Vert_\infty
+ \varepsilon^\gamma 
$$
and $\varphi_1(r):= \left( 2Kr-r^2 \right)^{1/\alpha_p}$ for $r\in
[0,K]$. Then $v$ is given by  
\beqn
\label{b14}
v:= K - \left( K^2 - u_{k,\varepsilon}^{\alpha_p} \right)^{1/2}
\eeqn
and satisfies 
\beqn
\label{b15}
\frac{\varepsilon^{\gamma\alpha_p}}{2K} \le v \le K - \left( K^2 -
M^{\alpha_p} \right)^{1/2} \le M^{\alpha_p/2} 
\eeqn
by (\ref{ex4b}). Thanks to the bounds (\ref{b15}), we can find $\mu$
large enough such that $\varphi_1$ enjoys the following properties: 
\beqn
\label{b16}
0 \ge \left( \frac{\varphi_1''}{\varphi_1'} \right)'(v) \ge -
\frac{C_1(\mu)}{v^2}\,, 
\eeqn
\beqn
\label{b17}
0 \le \frac{\varphi_1''}{\varphi_1'}(v) \le \frac{C_2(\mu)}{v}\,, 
\eeqn
\beqn
\label{b18}
\left( \frac{\varphi_1''}{\varphi_1'} \right)'(v) +
\frac{\alpha_p}{1-\alpha_p}\ \left( \frac{\varphi_1''}{\varphi_1'}
\right)^2(v) \le - \frac{1+\alpha_p}{2\alpha_p}\ \frac{1}{Kv}\,. 
\eeqn 
We then infer from (\ref{b16}) and (\ref{b17}) that 
$$
\mathcal{R}_{11} \ge - \frac{C_3(\mu)}{v^2}\ g^{p-4}\,.
$$
Therefore, by (\ref{b15}) and the elementary inequality $g\ge
\vert\nabla u_{k,\varepsilon}\vert$, we have 
$$
w\ \mathcal{R}_{11} \ge - \frac{\vert\nabla u_{k,\varepsilon}
\vert^2}{(\varphi_1')^2(v)}\ \frac{C_3(\mu)}{v^2}\ g^{p-4} \ge
-\frac{C_4(\mu)}{M\ v^{2/\alpha_p}}\ g^{p-2} \ge
-\frac{C_5(\mu)}{\varepsilon^{2 \gamma}}\ g^{p-2}\,. 
$$
Combining the previous inequality with (\ref{b10}) and (\ref{b18}), we
obtain 
$$
w^2\ \mathcal{R}_1 \ge \frac{C_6(\mu)\ M^{-\alpha_p/2}}{v}\ g^{p-2}\
w^2 - C_5(\mu)\ \varepsilon^{2(1-\gamma)}\ g^{p-2}\ w\,. 
$$
Now, we have $g\le \left\Vert \nabla u_{0,k} \right\Vert_\infty +
\varepsilon$ by (\ref{ex6}) and 
$$
g^2 \ge \vert\nabla u_{k,\varepsilon}\vert^2 = (\varphi_1')^2(v) \ge
C_7(\mu)\ M\ v^{2(1-\alpha_p)/\alpha_p}\ w  
$$
by (\ref{b15}). The previous lower bound for $w^2\ \mathcal{R}_1$ then
gives 
$$
w^2\ \mathcal{R}_1 \ge \frac{C_8(\mu)\
M^{(p-2-\alpha_p)/2}}{v^{((p-1)\alpha_p-(p-2))/\alpha_p}}\ w^{(p+2)/2}
- C_9(\mu,k)\ \varepsilon^{2(1-\gamma)}\ w\,. 
$$
Since $(p-1)\alpha_p\ge (p-2)$ and $v\le M^{\alpha_p/2}$ by
(\ref{b15}), we end up with 
\beqn
\label{b19}
w^2\ \mathcal{R}_1 \ge C_{10}(\mu)\ M^{(2(p-2)-p\alpha_p)/2}\
w^{(p+2)/2} - C_9(\mu,k)\ \varepsilon^{2(1-\gamma)}\ w\,. 
\eeqn
Next, since $q>1$ and $g\ge\varepsilon$, we infer from the
monotonicity of $\varphi_1$ and (\ref{b17}) that $\mathcal{R}_2\ge
0$. Recalling (\ref{b2}) and (\ref{b19}) we have shown that  
$$
\mathcal{L}_1 w := \partial_t w - \mathcal{A} w - \mathcal{V}\cdot
\nabla w + 2\ C_{10}(\mu)\ M^{(2(p-2)-p\alpha_p)/2}\ w^{(p+2)/2} - 2\
C_9(\mu,k)\ \varepsilon^{2(1-\gamma)}\ w \le 0  
$$
in $Q_\infty$. It is then straightforward to check that 
$$
S_1(t) := \left( \frac{1 + 2\ C_9(\mu,k)\ \varepsilon^{1/4}}{p\
C_{10}(\mu)} \right)^{2/p}\ M^{(p\alpha_p - 2(p-2))/p}\ t^{-2/p}
$$
satisfies $\mathcal{L}_1 S_1\ge 0$ in $\left( 0,\varepsilon^{-1/4}
\right)\times\RR^N$. The comparison principle then ensures that
$w(t,x)\le S_1(t)$ for $(t,x)\in \left( 0,\varepsilon^{-1/4}
\right)\times\RR^N$. The estimate (\ref{b12}) then readily follows
with the help of (\ref{b15}).  

\medskip

To prove (\ref{b13}) we take $\varphi_2(r):= \beta_{p,q}\
r^{1/\beta_{p,q}}$, so that $v=(u/\beta_{p,q})^{\beta_{p,q}}$
satisfies 
\beqn
\label{b20}
\frac{\varepsilon^{\gamma \beta_{p,q}}}{\beta_{p,q}^{\beta_{p,q}}} \le
v \le \frac{M^{\beta_{p,q}}}{\beta_{p,q}^{\beta_{p,q}}} \;\;\mbox{
with }\;\; M:= \Vert u_{0,k}\Vert_\infty + \varepsilon^\gamma\,, 
\eeqn
by (\ref{ex4b}). Concerning $\mathcal{R}_1$, the computations are much
simpler than in the previous case and it follows from the definition
of $\beta_{p,q}$ and (\ref{ex6}) that  
\bear
\nonumber
w^2\ \mathcal{R}_1 & \ge & C_{11}\ \frac{\beta_{p,q} -
\alpha_p}{\alpha_p \beta_{p,q}}\ \frac{g^{p-2}\ w^2}{v^2} - C_{12}\
\varepsilon^{(2\beta_{p,q}-\gamma)/\beta_{p,q}}\ g^{p-2}\ w \\ 
\label{b21}
w^2\ \mathcal{R}_1 & \ge & -C_{13}(k)\
\varepsilon^{(2\beta_{p,q}-\gamma)/\beta_{p,q}}\ w\,. 
\eear

For $\mathcal{R}_2$, we first claim that
\beqn
\label{b22}
(q-1)\ g^q + \varepsilon^q - q\ \varepsilon^2\ g^{q-2} \ge
(q-1-\varepsilon)\ g^q - C_{14}\ \left( \varepsilon^{(q+2)/2} +
\varepsilon^q \right)\,.
\eeqn
Indeed, if $q>2$, it follows from the Young inequality that
\bean
(q-1)\ g^q + \varepsilon^q - q\ \varepsilon^2\ g^{q-2} & \ge & (q-1)\ g^q
- \varepsilon\ g^q - 2\ (q-2)^{(q-2)/2}\ \varepsilon^{(q+2)/2} \\
& \ge & (q-1-\varepsilon)\ g^q - 2\ (q-2)^{(q-2)/2}\
\varepsilon^{(q+2)/2}\,. 
\eean
If $q\in (1,2]$, we have 
$$
(q-1)\ g^q + \varepsilon^q - q\ \varepsilon^2\ g^{q-2} \ge (q-1)\ g^q
+ \varepsilon^q - q\ \varepsilon^q \ge (q-1-\varepsilon)\ g^q - (q-1)\
\varepsilon^q\,,
$$
which completes the proof of (\ref{b22}). We then infer from
(\ref{b11}), (\ref{b20}), and (\ref{b22}) that 
\bean
\mathcal{R}_2 & \ge & \frac{1-\beta_{p,q}}{\beta_{p,q}}\
\frac{1}{v^{1/\beta_{p,q}}}\ \left[ (q-1-\varepsilon)\
(\varphi_2')^q(v)\ w^{q/2} - C_{14}\ \left( \varepsilon^{(q+2)/2} +
\varepsilon^q \right) \right]\\
& \ge & \frac{1-\beta_{p,q}}{\beta_{p,q}}\ (q-1-\varepsilon)\ 
v^{(q(1-\beta_{p,q})-1)/\beta_{p,q}}\ w^{q/2} - C_{15}\ \left(
\varepsilon^{(q+2-2\gamma)/2} + \varepsilon^{q-\gamma} \right) \\ 
& \ge & \frac{1-\beta_{p,q}}{\beta_{p,q}^{q(1-\beta_{p,q})}}\
(q-1-\varepsilon)\ M^{q(1-\beta_{p,q})-1}\ w^{q/2} - C_{15}\ \left(
\varepsilon^{(q+2-2\gamma)/2} + \varepsilon^{q-\gamma} \right)\,,
\eean
 
Recalling (\ref{b21}) we have thus shown that $w$ satisfies
\bean
\mathcal{L}_2 w & := & \partial_t w - \mathcal{A} w - \mathcal{V}\cdot
\nabla w + 2\ \frac{1-\beta_{p,q}}{\beta_{p,q}^{q(1-\beta_{p,q})}}\
(q-1-\varepsilon)\ M^{q(1-\beta_{p,q})-1}\ w^{(q+2)/2} \\
& - & C_{16}(k)\ \omega(\varepsilon)\ w \le 0
\eean
in $Q_\infty$, where
$\omega(\varepsilon):=\varepsilon^{(2\beta_{p,q}-\gamma)/\beta_{p,q}} +
\varepsilon^{(q+2-2\gamma)/2} + \varepsilon^{q-\gamma}\to 0$ as
$\varepsilon\to 0$ by the choice (\ref{ex3b}) of $\gamma$. 
The function  
$$
S_2(t) := \frac{\beta_{p,q}^{2(1-\beta_{p,q})}}{2^{2/q}\ (1-\beta_{p,q})^{2/q}\
(q-1-\varepsilon)^{2/q}}\ \left( \frac{2 + q\ C_{16}(k)\
\omega(\varepsilon)^{1/2}}{q} \right)^{2/q}\
M^{2(1-q(1-\beta_{p,q}))/q}\ t^{-2/q} 
$$
satisfies $\mathcal{L}_2 S_2\ge 0$ in $\left( 0,\omega(\varepsilon)^{-1/2}
\right)\times\RR^N$. We then deduce from the comparison principle that
$w(t,x)\le S_2(t)$ for $(t,x)\in \left( 0,\omega(\varepsilon)^{-1/2}
\right)\times\RR^N$. The estimate (\ref{b13}) then readily follows. \qed


\section{Existence}\label{exist}
\setcounter{thm}{0}
\setcounter{equation}{0}

We are now in a position to prove Theorem~\ref{tha1} and proceed along the
lines of \cite{GGK03}. 

\noindent\textbf{Step~1: $\varepsilon\to 0$}. We first let
$\varepsilon\to 0$. For that purpose, we observe that the gradient
bound (\ref{ex6}) and (\ref{ex4}) imply the time equicontinuity of
$(u_{k,\varepsilon})_{\varepsilon>0}$.

\begin{lem}\label{lec1}
For $k\ge 1$, $\varepsilon>0$, $x\in\RR^N$, $t_1\ge 0$, and $t_2>t_1$,
we have 
$$
\left| u_{k,\varepsilon}(t_2,x) - u_{k,\varepsilon}(t_1,x) \right| \le
C\ \left( \|\nabla u_{0,k}\|_\infty + \|\nabla u_{0,k}\|_\infty^{p-1}
\right)\ (t_2-t_1)^{1/2} + \|\nabla u_{0,k}\|_\infty^q\ (t_2-t_1)\,.
$$
\end{lem}

The proof of Lemma~\ref{lec1} is similar to that of
\cite[Lemma~5]{GGK03} to which we refer.

\medskip

We next fix $k\ge 1$. Owing to (\ref{ex4b}), (\ref{ex6}), and
Lemma~\ref{lec1}, we may apply the Arzel\`a-Ascoli theorem to obtain a
subsequence of $(u_{k,\varepsilon})_{\varepsilon>0}$ (not relabeled)
and a non-negative function
$u_k\in\mathcal{BC}([0,\infty)\times\RR^N)$ such that 
\beqn
\label{c1}
u_{k,\varepsilon} \longrightarrow u_k \;\;\mbox{ uniformly on any
compact subset of }\;\; [0,\infty)\times\RR^N\,.
\eeqn
Furthermore, as $u_{k,\varepsilon}$ is a classical solution to
(\ref{ex4}), (\ref{ex5}), the classical stability result for
continuous viscosity solutions allows us to conclude that $u_k$ is a
viscosity solution to (\ref{vhj1}) with initial condition $u_{0,k}$
(see, e.g., \cite[Theorem~1.4]{CEL84} or
\cite[Th\'eor\`eme~2.3]{Bl94}). By (\ref{c1}) and weak convergence
arguments, we next infer from (\ref{ex4b}), (\ref{b12}), and
(\ref{b13}) that  
\bear
\label{c1b}
0 \le u_k(t,x) & \le & \|u_0\|_\infty\,, \\
\label{c2}
\left\vert \nabla\left( u_k^{\alpha_p} \right)(t,x) \right\vert & \le
& C\ \Vert u_{0,k}\Vert_\infty^{(p \alpha_p+2-p)/p}\ t^{-1/p}\,, \\ 
\nonumber
 & & \\
\label{c3}
\left\vert \nabla\left( u_k^{\beta_{p,q}} \right)(t,x) \right\vert &
\le & \frac{\beta_{p,q}}{(q^2-q)^{1/q} (1-\beta_{p,q})^{1/q}}\ \Vert
u_{0,k}\Vert_\infty^{(q \beta_{p,q}+1-q)/q}\ t^{-1/q} 
\eear
for all $(t,x)\in Q_\infty$. Finally, (\ref{ex4}) also reads 
$$
\partial_t u_{k,\varepsilon} - \mbox{ div } \left( \vert\nabla
u_{k,\varepsilon} \vert^{p-2}\ \nabla u_{k,\varepsilon} \right)  =
\mbox{ div } \left( f_{k,\varepsilon} \right) + g_{k,\varepsilon}
\;\;\mbox{ in }\;\; Q_\infty 
$$
with
$$
f_{k,\varepsilon}:= \left\{ a_\varepsilon\left( \vert\nabla
u_{k,\varepsilon} \vert^2\right) - \vert\nabla u_{k,\varepsilon}
\vert^{p-2} \right\}\ \nabla u_{k,\varepsilon} \;\;\mbox{ and }\;\;
g_{k,\varepsilon}:= - b_\varepsilon\left( \vert\nabla
u_{k,\varepsilon}\vert^2 \right)\,. 
$$
It follows from the definition of $a_\varepsilon$ and (\ref{ex6}) that
$(g_{k,\varepsilon})$ is bounded in $L^\infty(Q_\infty)$ and
$(f_{k,\varepsilon})$ converges to zero in $L^\infty(Q_\infty)$ as
$\varepsilon\to 0$. We may then apply \cite[Theorem~4.1]{BM92} to
conclude that 
\beqn
\label{c4}
\nabla u_{k,\varepsilon} \longrightarrow \nabla u_k \;\;\mbox{ a.e. in
}\;\; Q_\infty\,. 
\eeqn
Consequently, upon extracting a further subsequence, we may assume that 
\beqn
\label{c5}
\nabla u_{k,\varepsilon} \longrightarrow \nabla u_k \;\;\mbox{ a.e. in
}\;\; L^r((0,T)\times B(0,R)) 
\eeqn
for every $r\in [1,\infty)$, $T>0$, and $R>0$. It then readily follows
that $u_k$ satisfies (\ref{weakf}) with $u_{0,k}$ instead of $u_0$.  

\medskip

\noindent\textbf{Step~2: $k\to\infty$}. It remains to pass to the
limit as $k\to\infty$. To this end we first observe that (\ref{ex1})
implies that $u_{0,k}(x)-u_{0,k+1}(y) \le \|\nabla u_{0,k}\|_\infty\
|y-x|$ for $k\ge 1$, $x\in\RR^N$, and $y\in\RR^N$. It then follows
from the comparison principle \cite[Theorem~2.1]{GGIS91} that  
\beqn
\label{c6}
u_k(t,x) \le u_{k+1}(t,x) \;\;\mbox{ for } (t,x)\in Q_\infty
\;\;\mbox{ and }\;\; k\ge 1\,. 
\eeqn
Therefore, by (\ref{ex1}), (\ref{c1b}), and (\ref{c6}), the function
\beqn
\label{c7}
u(t,x) := \sup_{k\ge 1} u_k(t,x) \in \left[ 0 , \|u_0\|_\infty \right]
\eeqn
is well-defined for $(t,x)\in [0,\infty)\times\RR^N$. We next readily
deduce from (\ref{c1b}) and (\ref{c2}) that, for $\tau>0$,  
\beqn
\label{c8}
\|\nabla u_k(t)\|_\infty \le C\ \|u_0\|_\infty^{2/p}\ t^{-1/p} \le C\
\|u_0\|_\infty^{2/p}\ \tau^{-1/p} \;\;\mbox{ for }\;\; t\ge \tau\,. 
\eeqn
Thanks to (\ref{c8}) we may argue as in the previous step and conclude that 
\beqn
\label{c9}
u_k \longrightarrow u \;\;\mbox{ uniformly on any compact subset of
}\;\; Q_\infty\,. 
\eeqn
Using again the stability of continuous viscosity solutions, we deduce
from the convergence (\ref{c9}) that $(t,x)\longmapsto u(t+\tau,x)$ is
a viscosity solution to
(\ref{vhj1}) with initial condition $u(\tau)$ for each $\tau>0$. In
addition, denoting by $\tilde{u}_k$ the solution to the $p$-Laplacian
equation (\ref{ple}) with initial condition $u_{0,k}$, the comparison
principle entails that 
\beqn
\label{c10}
u_k(t,x) \le \tilde{u}_k(t,x) \;\;\mbox{ for }\;\; (t,x)\in Q_\infty
\;\;\mbox{ and }\;\; k\ge 1\,. 
\eeqn
Furthermore, $(\tilde{u}_k)_{k\ge 1}$ converges uniformly on any
compact subset of $[0,\infty)\times\RR^N$ towards the solution
$\tilde{u}$ to the $p$-Laplacian equation (\ref{ple}) with initial
condition $u_0$ \cite[Ch.~III]{DB93}. This property and (\ref{c10})
warrant that $u(t,x)\le\tilde{u}(t,x)$ for $(t,x)\in
[0,\infty)\times\RR^N$. Recalling (\ref{c7}), we thus obtain the
following inequality 
\beqn
\label{c11}
u_k(t,x)\le u(t,x) \le \tilde{u}(t,x) \;\;\mbox{ for }\;\; (t,x)\in
Q_\infty \;\;\mbox{ and }\;\; k\ge 1\,. 
\eeqn
We then infer from (\ref{c11}) that $(u(.+1/j))_{j\ge 1}$ converges
towards $u$ uniformly on any compact subset of $[0,\infty)\times\RR^N$
as $j\to\infty$. Using once more the stability of continuous viscosity
solutions, we conclude that $u$ is a viscosity solution to
(\ref{vhj1}), (\ref{vhj2}). We next argue as in the previous step to
deduce from (\ref{c2}) and (\ref{c3}) that $u$ satisfies (\ref{ge1}),
(\ref{ge2}) and (\ref{weakf}) for $t>s>0$. In addition, $u\in
L^\infty(Q_\infty)$ by (\ref{se1}) and we deduce from (\ref{se1}) and
(\ref{ge1}) that $\|\nabla u(t)\|_\infty \le C\ \|u_0\|_\infty^{2/p}\
t^{-1/p}$ for $t>0$.  Consequently, $\nabla u$ belongs to
$L^{p-1}((0,T)\times B(0,R))$ for all $T>0$ and $R>0$. We then let
$s\to 0$ in (\ref{weakf}) to conclude that $\nabla u \in
L^q((0,T)\times B(0,R))$ for all $T>0$ and $R>0$ which in turn
warrants that (\ref{weakf}) is also valid for $s=0$.  

To complete the proof of Theorem~\ref{tha1}, it remains to check the
uniqueness assertion for $u_0\in\mathcal{BUC}(\RR^N)$ which actually
follows at once from \cite[Theorem~2.1]{GGIS91}.  


\section{Temporal decay estimates}\label{tde}

\setcounter{thm}{0}
\setcounter{equation}{0}

This section is devoted to the proof of Proposition~\ref{pra3}. Let
us start with the following lemma:
\begin{lem}
Let $u$ be a solution of (\ref{vhj1}), (\ref{vhj2}). If $t>s\geq 0$, then
\bear
\label{decpl}
\Vert \nabla u(t)\Vert_\infty & \le & C\ \Vert u(s)\Vert_\infty^{2/p}\
(t-s)^{-1/p}\,,\\ 
\label{dechj}
\Vert \nabla u(t)\Vert_\infty & \le & C \Vert u(s)\Vert_\infty^{1/q}\
(t-s)^{-1/q}\, . 
\eear
\end{lem}

\noindent\textbf{Proof.} We write
\[
\vert \nabla u(t) \vert=\frac{1}{\gamma}\ u^{1-\gamma}\ \vert \nabla u^\gamma\vert 
\]
for $\gamma=\alpha_p$ and $\gamma=\beta_{p,q}$ and use the estimates
(\ref{ge1}) and (\ref{ge2}).
\qed

\medskip

\noindent\textbf{Proof of Proposition~\ref{pra3}.} We first prove
(\ref{dec1}). Combining the Gagliardo-Nirenberg inequality, the time
monotonicity of $\Vert u\Vert_1$ and the previous lemma, we obtain 
\begin{eqnarray*}
\Vert u(t)\Vert_\infty^q&\leq& C\ \Vert\nabla
u(t)\Vert_\infty^{qN/(N+1)}\ \Vert u(t)\Vert_1^{q/(N+1)}\\ 
&\leq& C\Vert\nabla u(t)\Vert_\infty^{qN/(N+1)}\ \Vert u_0\Vert_{1}^{q/(N+1)}\\
&\leq& C(t-s)^{-N/(N+1)}\ \Vert u(s)\Vert_\infty^{N/(N+1)}\ \Vert
u_0\Vert_{1}^{q/(N+1)}\,. 
\end{eqnarray*}
Integrating with respect to $t$ over $(s,\infty)$, we obtain
\begin{eqnarray*}
\tau(s):= \int_s^\infty \frac{\Vert u(t)\Vert_{\infty}^q}{t}\ dt
&\leq& C\ \Vert u(s)\Vert_{\infty}^{N/(N+1)}\ \Vert
u_0\Vert_{1}^{q/(N+1)}\ \int_s^\infty \frac{dt}{(t-s)^{N/(N+1)} t}\\ 
&\leq& C\ s^{-N/(N+1)}\ \Vert u_0\Vert_{1}^{q/(N+1)}\ \Vert
u(s)\Vert_{\infty}^{N/(N+1)}\,, 
\end{eqnarray*}
whence
\[
\tau(s) \leq C\ \Vert u_0\Vert_{1}^{q/(N+1)}\, \left( -\tau'(s)
\right)^{N/q(N+1)}\ s^{-(N(q-1))/q(N+1)}\,. 
\]
Introducing $\tilde\tau(s)=\tau(s^{1/q})$ gives
\[
\frac{d\tilde\tau}{ds}(s) + C\ \Vert u_0\Vert_{1}^{-q^2/N}\
\tilde\tau(s)^{q(N+1)/N}\leq 0\,. 
\]
A direct computation shows that $\tilde\tau(s)\leq C\ \Vert
u_0\Vert_1^{q^2\xi}\ s^{-N\xi}$ from which we deduce that 
\[
 \tau(s)\leq C\ \Vert u_0\Vert_1^{q^2\xi}\ s^{-qN\xi}\,.
\]
Now, using the time monotonicity of $\Vert u\Vert_\infty$, we obtain
\[
C\ s^{-qN\xi}\ \Vert u_0\Vert_1^{q^2\xi} \geq \tau(s) \geq
\int_s^{2s}\frac{\Vert u(t)\Vert_{\infty}^q}{t}\,dt \geq
\int_s^{2s}\frac{\Vert u(2s)\Vert_{\infty}^q}{t}\,dt=\ln(2)\ \Vert
u(2s)\Vert_{\infty}^q\,, 
\]
whence (\ref{dec1}). The estimate (\ref{dec2}) then readily follows
from (\ref{dec1}) by (\ref{dechj}). A similar proof relying on
(\ref{decpl}) gives the estimates (\ref{dec3}) and (\ref{dec4}). \qed 


\section{Limit values of $\Vert u(t)\Vert_1$}\label{lvl1}

\setcounter{thm}{0}
\setcounter{equation}{0}

In this section we investigate the possible values of the limit as
$t\to\infty$ of the $L^1$-norm of non-negative solutions to
(\ref{vhj1}), (\ref{vhj2}) and prove Proposition~\ref{pra4}. We first
show that, if $q$ is small enough, the dissipation mechanism induced
by the nonlinear absorption term is sufficiently strong to drive the
$L^1$-norm of $u$ to zero in infinite time.  

\begin{prop}\label{pre1}
If $q\in (1,q_*]$ then
$$
\lim_{t\to\infty} \Vert u(t)\Vert_1=0\,.
$$
\end{prop}

\noindent\textbf{Proof.} It first follows from the integration 
of (\ref{vhj1}) over $(0,t)\times\RR^N$ that  
\beqn
\label{e1}
\Vert u(t)\Vert_1 + \int_0^t \Vert\nabla u(s)\Vert_q^q\ ds = 
\Vert u_0\Vert_1\,,
\eeqn
which readily implies that $t\longmapsto \Vert\nabla u(t)\Vert_q^q$ 
belongs to $L^1(0,\infty)$. Consequently, 
\beqn
\label{e2}
\omega(t):= \int_t^\infty \Vert\nabla u(s)\Vert_q^q\ ds
\mathop{\longrightarrow}_{t\to\infty} 0\,. 
\eeqn

We next consider a $\mathcal{C}^\infty$-smooth function $\vartheta$ in
$\RR^N$ such that $0\le \vartheta \le 1$ and 
$$
\vartheta(x) = 0 \;\;\mbox {if }\;\; \vert x\vert\le 1/2 \;\;\mbox{
and }\;\; \vartheta(x) = 1 \;\;\mbox{ if }\;\; \vert x\vert\ge 1\,.
$$
For $R>0$ and $x\in\RR^N$ we put $\vartheta_R(x) = \vartheta(x/R)$. We
multiply (\ref{vhj1}) by $\vartheta_R(x)$ and integrate over
$(t_1,t_2)\times\RR^N$ to obtain 
$$
\int_{\RR^N} u(t_2,x)\ \vartheta_R(x)\ dx \le \int_{\RR^N} u(t_1,x)\
\vartheta_R(x)\ dx + \frac{1}{R}\ \int_{t_1}^{t_2} \vert\nabla
u(s,x)\vert^{p-2}\ \nabla\vartheta\left( \frac{x}{R} \right)\ \nabla
u(s,x)\ dxds\,,  
$$
which, together with the properties of $\vartheta$, gives 
\beqn
\label{e3}
\int_{\{ \vert x\vert\ge 2R\}} u(t_2,x)\ dx \le \int_{\{ \vert
x\vert\ge R\}} u(t_1,x)\ dx + \frac{1}{R}\ \int_{t_1}^{t_2} \int_{\RR^N} 
\left\vert\nabla\vartheta\left( \frac{x}{R} \right)\right\vert\ 
\vert\nabla u(s,x)\vert^{p-1}\ dxds\,.
\eeqn

\medskip

\noindent\textsl{Case~1: $q\in [p-1,q_*]$}. By the H\"older inequality
we have  
\bean
& & \frac{1}{R}\ \int_{t_1}^{t_2} \int_{\RR^N} \left\vert\nabla\vartheta\left(
\frac{x}{R} \right)\right\vert\ \vert\nabla u(s,x)\vert^{p-1}\ dxds \\
& \le & R^{(N(q-p+1)-q)/q}\ (t_2-t_1)^{(q-p+1)/q}\
\Vert\nabla\vartheta\Vert_{(q-p+1)/q}\ \left( \int_{t_1}^{t_2}
\Vert\nabla u(s)\Vert_q^q\ dxds \right)^{(p-1)/q} \\ 
& \le & C\ R^{(N(q-p+1)-q)/q}\ \omega(t_1)^{(p-1)/q}\
(t_2-t_1)^{(q-p+1)/q}\,. 
\eean
Combining the above inequality with (\ref{dec1}), (\ref{e3}) and the
time monotonicity of $\Vert u\Vert_1$ we obtain  
\bean
\Vert u(t_2)\Vert_1 & = & \int_{\{ \vert x\vert\le 2R\}} u(t_2,x)\ dx
+ \int_{\{ \vert x\vert\ge 2R\}} u(t_2,x)\ dx \\ 
& \le & C\ R^N\ \Vert u(t_2)\Vert_\infty + \int_{\{ \vert x\vert\ge
R\}} u(t_1,x)\ dx \\
& + & C\ R^{(N(q-p+1)-q)/q}\ \omega(t_1)^{(p-1)/q}\
(t_2-t_1)^{(q-p+1)/q} \\ 
& \le & \int_{\{ \vert x\vert\ge R\}} u(t_1,x)\ dx + C\ R^N\
(t_2-t_1)^{-N\xi} \\
& + & C\ R^{(N(q-p+1)-q)/q}\ \omega(t_1)^{(p-1)/q}\
(t_2-t_1)^{(q-p+1)/q}\,. 
\eean
Choosing
$$
R = R(t_1,t_2) := \omega(t_1)^{(p-1)/(q+N(p-1))}\
(t_2-t_1)^{(qN\xi+q-p+1)/(q+N(p-1))} 
$$
we are led to 
\bean
\Vert u(t_2)\Vert_1 & \le & \int_{\{ \vert x\vert\ge R(t_1,t_2)\}}
u(t_1,x)\ dx \\ 
& + & C\ \omega(t_1)^{(N(p-1))/(q+N(p-1))}\
(t_2-t_1)^{-qN\xi(N+1)(q_*-q)/(q+N(p-1))}\,. 
\eean
Since $\xi>0$ and $q_*-q>0$ we may let $t_2\to\infty$ in the previous
inequality to conclude that 
\bean
I_1(\infty) & \le & 0 \;\;\mbox{ if }\;\; q\in [p-1,q_*)\,, \\
I_1(\infty) & \le & C\ \omega(t_1)^{(N(p-1))/(q_*+N(p-1))} \;\;\mbox{
if }\;\; q=q_*\,. 
\eean
We have used here that $R(t_1,t_2)\to\infty$ as $t_2\to\infty$ and
that $u(t_1)\in 
L^1(\RR^N)$. Owing to the non-negativity of $I_1(\infty)$, we readily
obtain that $I_1(\infty)=0$ if $q\in [p-1,q_*)$. When $q=q_*$,
we let $t_1\to\infty$ and use (\ref{e2}) to conclude that
$I_1(\infty)=0$ also in that case. 

\noindent\textsl{Case~2: $q\in (1,p-1)$}. By (\ref{dec2}) and
(\ref{e3}) we have  
\bean
\int_{\{ \vert x\vert\ge 2R\}} u(t_2,x)\ dx & \le & \int_{\{ \vert
x\vert\ge R\}} u(t_1,x)\ dx + \frac{1}{R}\
\Vert\nabla\vartheta\Vert_\infty\ \int_{t_1}^{t_2} \Vert\nabla
u(s)\Vert_\infty^{p-1-q}\ \Vert\nabla u(s)\Vert_q^q\ ds \\ 
& \le & \int_{\{ \vert
x\vert\ge R\}} u(t_1,x)\ dx + \frac{C}{R}\ \int_{t_1}^{t_2}
s^{-(p-1-q)(N+1)\xi}\ \Vert\nabla u(s)\Vert_q^q\ ds \\ 
& \le & \int_{\{ \vert
x\vert\ge R\}} u(t_1,x)\ dx + \frac{C}{R}\ t_1^{-(p-1-q)(N+1)\xi}\
\omega(t_1)\,. 
\eean
Taking $t_1=1$ and noting that $\omega(t_1)\le\omega(0)\le \Vert
u_0\Vert_1$, we end up with 
$$
\int_{\{ \vert x\vert\ge 2R\}} u(t_2,x)\ dx \le \int_{\{ \vert
x\vert\ge R\}} u(1,x)\ dx + \frac{C}{R}\,, \quad t_2\ge 1\,.
$$
We then infer from (\ref{dec1}) and the above inequality that, if
$t_2\ge 1$,  
$$
\Vert u(t_2)\Vert_1 \le C\ R^N\ t^{-N\xi} + \int_{\{ \vert
x\vert\ge R\}} u(1,x)\ dx + \frac{C}{R}
$$
and the choice $R=R(t_2)=t_2^{(N\xi)/(N+1)}$ gives 
$$
\Vert u(t_2)\Vert_1 \le \int_{\{ \vert
x\vert\ge R(t_2)\}} u(1,x)\ dx + C\ t_2^{-(N\xi)/(N+1)}\,.
$$
Since $R(t_2)\to\infty$ as $t_2\to\infty$ and $u(1)\in L^1(\RR^N)$ we
may let $t_2\to\infty$ in the above inequality to establish that
$I_1(\infty)=0$, which completes the proof of
Proposition~\ref{pre1}. \qed  

\medskip

We next turn to higher values of $q$ and adapt an argument of
\cite[Theorem~6]{BL99} to show the positivity of $I_1(\infty)$.  

\begin{prop}\label{pre2}
Assume that $\Vert u_0\Vert_1>0$ and $q>q_*$. Then $I_1(\infty)>0$. 
\end{prop}

\noindent\textbf{Proof.} Since $u_0\in\mathcal{BC}(\RR^N)$ is not
identically equal to zero there are $x_0\in\RR^N$ and a radially
symmetric and non-increasing continuous function $U_0\not\equiv 0$
such that $u_0(x)\ge U_0(x-x_0)$. Denoting by $U$ the solution to
(\ref{vhj1}) with initial condition $U_0$ it follows from the
invariance of (\ref{vhj1}) by translation and the comparison principle
that  
\beqn
\label{e4}
u(t,x) \ge U(t,x-x_0) \,, \quad (t,x)\in [0,\infty)\times\RR^N\,.
\eeqn

\smallskip

Let $\tau>0$ and $x\in\RR^N$. 
Since 
$$
\nabla U(\tau,x) = \frac{p-1}{p-2}\ U(\tau,x)^{1/(p-1)}\ \nabla\left(
U^{(p-2)/(p-1)} \right)(\tau,x) 
$$ 
and $q>q_*>p-1$, we infer from (\ref{ge1r}) and the time monotonicity
of $\Vert u\Vert_\infty$ that 
\bean
\vert\nabla U(\tau,x)\vert^q & \le & \left( \frac{p-1}{p-2} \right)^q\
U(\tau,x)^{q/(p-1)}\ \left\vert \nabla\left(  U^{(p-2)/(p-1)}
\right)(\tau,x) \right\vert^q \\ 
& \le & C\ U(\tau,x)\ \Vert U(\tau)\Vert_\infty^{(q-p+1)/(p-1)}\
\left\Vert U\left( \frac{\tau}{2} \right) \right\Vert_
\infty^{q(p-2)/p(p-1)}\ \tau^{-q/p} \\ 
& \le & C\ U(\tau,x)\ \left\Vert U\left( \frac{\tau}{2} \right)
\right\Vert_ \infty^{(2q-p)/p}\ \tau^{-q/p}\,, 
\eean
whence, by (\ref{dec3}), 
\beqn
\label{e5}
\vert\nabla U(\tau,x)\vert^q \le C\ U(\tau,x)\ \tau^{-\eta/\xi}\,.
\eeqn

Consider now $s\in (0,\infty)$ and $t\in (s,\infty)$. It follows
from (\ref{vhj1}) and (\ref{e5}) that 
\bean
\Vert U(t)\Vert_1 & = & \Vert U(s)\Vert_1 
- \int_s^t \int_{\RR^N} \vert\nabla U(\tau,x)\vert^q\ dxd\tau
\\ 
& \ge & \Vert U(s)\Vert_1 - C\ \int_s^t \tau^{-\eta/\xi}\ \Vert
U(\tau)\Vert_1\ d\tau \,. 
\eean
Owing to the monotonicity of $\tau\longmapsto \Vert U(\tau)\Vert_1$,
we further obtain 
$$
\Vert U(t)\Vert_1 \ge \Vert U(s)\Vert_1 \ \left({ 1 - C\ \int_s^t
\tau^{-\eta/\xi}\ d\tau }\right)\,. 
$$
Since $q>q_*$ we have $\eta>\xi$ and the right-hand side of the above
inequality has a finite limit as $t\to\infty$. We may then let $t\to
\infty$ to obtain  
$$
\mathcal{I}_1(\infty):= \lim_{t\to\infty} \Vert U(t)\Vert_1 \ge \Vert
U(s)\Vert_1\ \left({ 1 - C\ s^{-(\eta-\xi)/\xi} }\right)\,, \quad s>0\,. 
$$
Consequently, for $s$ large enough, we have $\mathcal{I}_1(\infty)\ge
\Vert U(s)\Vert_1/2$, while \cite[Lemma~4.1]{ATU04} warrants that
$\Vert U(s)\Vert_1>0$ for each $s\ge 0$ since $U_0\not\equiv
0$. Therefore, $\mathcal{I}_1(\infty)>0$. Recalling (\ref{e4}) we
realize that $\Vert u(t)\Vert_1\ge \Vert U(t)\Vert_1$ for each $t\ge
0$ so that $I_1(\infty)\ge \mathcal{I}_1(\infty) > 0$. \qed


\section{Compactly supported initial data}\label{csid}

\setcounter{thm}{0}
\setcounter{equation}{0}

This section is devoted to the proofs of Theorem~\ref{tha5} and
Corollary~\ref{cora6}. Let $u_0\in L^1(\RR^N)\cap\mathcal{BC}(\RR^N)$
be a non-negative initial condition with compact support in the ball
$B(0,R_0)$ for some $R_0>0$. Denoting by $u$ the corresponding
solution to (\ref{vhj1}), (\ref{vhj2}) and by $v$ the corresponding 
solution to the $p$-Laplacian equation  
\beqn
\label{ple1}
\partial_t v - \Delta_p v = 0\ \,, \quad (t,x)\in Q_\infty\,,
\eeqn
with initial condition $v(0)=u_0$, the comparison principle ensures that 
\beqn
\label{f1}
0 \le u(t,x) \le v(t,x)\,, \quad (t,x)\in Q_\infty\,.
\eeqn
Since $u_0$ is compactly supported, so is $v(t)$ for each $t\ge 0$ by
\cite[Lemma~8.1]{DB93} and $\mbox{ Supp } v(t) \subset B(0, C_1
t^\eta)$. Consequently, $u(t)$ is compactly supported for each $t\ge
0$ with $\mbox{ Supp } u(t) \subset B(0, C_1 t^\eta)$. In particular,
the support of $u$ does not expand faster than that of $v$ with
time. A natural question is then whether the damping term slows down
this expansion and the answer depends heavily on the value of $q$. We
shall thus distinguish between three cases in the proof of
Theorem~\ref{tha5}. 

We first note that, since $u_0$ is non-negative continuous and
compactly supported, there exists a non-negative continuous radially
symmetric and non-increasing function $U_0$ with compact support such
that $0\le u_0\le U_0$. Denoting by $U$ the corresponding solution to
(\ref{vhj1}) with initial condition $U(0)=U_0$, the function
$x\longmapsto U(t,x)$ is also radially symmetric and non-increasing
for each $t\ge 0$ and we deduce from the comparison principle that  
\beqn
\label{f2}
0 \le u(t,x) \le U(t,x)\,, \quad (t,x)\in [0,\infty)\times\RR^N\,.
\eeqn
Moreover, by comparison with the $p$-Laplacian equation, $U(t)$ is
also compactly supported for each $t\ge 0$ with $\mbox{ Supp } U(t)
\subset B(0, \sigma(t))$ for some $\sigma(t)>0$. Clearly,  
\beqn
\label{f3}
\varrho(t) \le \sigma(t)\,, \quad t\ge 0\,,
\eeqn
by (\ref{f2}). 

It next follows from (\ref{vhj1}) that, if $y$ is a non-negative 
function in $\mathcal{C}^1([0,\infty))$, we have 
\bean
\frac{d}{dt} \int_{\{ \vert x\vert\ge y(t)\}} U(t,x)\ dx & = &
\int_{\{ \vert x\vert\ge y(t)\}} \partial_t U(t,x)\ dx - y'(t)\
\int_{\{ \vert x\vert= y(t)\}} U(t,x)\ dx \\ 
& \le & \int_{\{ \vert x\vert\ge y(t)\}} \mbox{ div } \left(
\vert\nabla U\vert^{p-2}\ \nabla U \right)(t,x)\ dx \\ 
& - & y'(t)\ \int_{\{ \vert x\vert= y(t)\}} U(t,x)\ dx \\
& \le & - \int_{\{ \vert x\vert= y(t)\}} \vert\nabla U(t,x)\vert^{p-2}\
\nabla U(t,x)\cdot \frac{x}{\vert x\vert}\ dx \\ 
& - & y'(t)\ \int_{\{ \vert x\vert= y(t)\}} U(t,x)\ dx \\ 
& \le & \int_{\{ \vert x\vert= y(t)\}} \left\{ \left\vert\nabla
U(t,x)\right\vert^{p-1} - y'(t)\ U(t,x) \right\}\ dx\,, 
\eean
\bear
\nonumber
& & \frac{d}{dt} \int_{\{ \vert x\vert\ge y(t)\}} U(t,x)\ dx \\
\label{f4}
& \le & \int_{\{ \vert x\vert=y(t)\}} \left\{ \frac{p-1}{p-2}\
\left\vert\nabla \left( U^{(p-2)/(p-1)} \right)(t,x)\right\vert^{p-1}
- y'(t) \right\}\ U(t,x)\ dx\,. 
\eear
The next step is to use the gradient estimates established in
Theorem~\ref{tha2} to find a suitable function $y$ for which the
right-hand side of (\ref{f4}) is non-positive. The gradient estimates
depending on the value of $q$, we handle separately the cases $q\in
(1,p-1]$ and $q\in (p-1,q_*)$. 

\medskip

\noindent\textbf{Proof of Theorem~\ref{tha5}: $q\in (1,p-1]$.} In that
case we infer from (\ref{ge2rb}) and (\ref{dec1}) that 
\bean
\left\vert\nabla \left( U^{(p-2)/(p-1)} \right)(t,x)\right\vert^{p-1}
& \le & C\ \left\Vert u\left( \frac{t}{2} \right)
\right\Vert_\infty^{(p-1-q)/q}\ t^{-(p-1)/q} \\ 
& \le & C\ t^{-\xi((p-1)(N+1)-N)}\,,
\eean
so that (\ref{f4}) becomes  
$$
\frac{d}{dt} \int_{\{ \vert x\vert\ge y(t)\}} U(t,x)\ dx \le \int_{\{
\vert x\vert=y(t)\}} \left\{ C\ t^{-\xi((p-1)(N+1)-N)} - y'(t)
\right\}\ U(t,x)\ dx 
$$
Choosing $y'(t):= C\ t^{-\xi((p-1)(N+1)-N)}$ for $t\ge 1$ and
$y(1)=\sigma(1)$, we conclude that 
$$
\int_{\{ \vert x\vert\ge y(t)\}} U(t,x)\ dx \le \int_{\{ \vert
x\vert\ge \sigma(1)\}} U(1,x)\ dx = 0 
$$
for $t\ge 1$. Consequently, $\sigma(t)\le y(t)$ for $t\ge 1$ from
which we deduce that $\varrho(t)\le y(t)$ for $t\ge 1$ by
(\ref{f2}). Now, either $q\in (1,p-1)$ and
$\xi((p-1)(N+1)-N)>1$. Therefore $y(t)$ has a finite limit as
$t\to\infty$ from which (\ref{rad1}) readily follows. Or $q=p-1$ and
$y(t)=\sigma(1) + C\ \ln{t}$ which gives (\ref{rad2}). \qed 

\medskip

We next consider the case $q\in (p-1,q_*)$ which turns out to be more
complicated as (\ref{ge2rb}) is no longer available. We instead use
(\ref{ge1r}) which somehow provides less information and thus
complicates the proof. We shall also need the following lemma which is
an easy consequence of the Poincar\'e and H\"older inequalities.  

\begin{lem}\label{lef1}
There is a positive constant $\kappa$ depending only on $N$ and $q$
such that, if $R>0$ and $w$ is a function in $W_0^{1,q}(B(0,R))$ then

\beqn
\label{f5}
R^{-1/\xi}\ \Vert w\Vert_{L^1(B(0,R))}^q \le \kappa\ \Vert\nabla
w\Vert_{L^q(B(0,R))}^q\,. 
\eeqn
\end{lem}

\medskip

\noindent\textbf{Proof of Theorem~\ref{tha5}: $q\in (p-1,q_*)$.} We
fix $t_0\ge 0$. It follows from (\ref{ge1r}) and (\ref{dec1}) that 
\bean
\frac{p-1}{p-2}\ \left\vert\nabla \left( U^{(p-2)/(p-1)}
\right)(t,x)\right\vert^{p-1} & \le & C\ \left\Vert u\left(
\frac{t+t_0}{2} \right) \right\Vert_\infty^{(p-2)/p}\
(t-t_0)^{-(p-1)/p} \\  
& \le & C\ \Vert u(t_0)\Vert_1^{q\xi(p-2)/p}\ (t-t_0)^{-(p-1+N\xi(p-2))/p}
\eean
for $t\ge t_0$. Since $q>p-1> N(p-1)/(N+1)$, we have $1-N\xi(p-2)>0$
and we choose $y(t) = \sigma(t_0) + p C\ \Vert u(t_0)\Vert_1^{q\xi(p-2)/p}\
(t-t_0)^{(1-N\xi(p-2))/p} / (1-N\xi(p-2))$ for $t\ge t_0$. The
previous inequality then reads
$$
\frac{p-1}{p-2}\ \left\vert\nabla \left( U^{(p-2)/(p-1)}
\right)(t,x)\right\vert^{p-1} \le y'(t)\,, \quad t\ge t_0\,.
$$
Combining the latter estimate with (\ref{f4}) we realize that 
$$
\frac{d}{dt} \int_{\{ \vert x\vert\ge y(t)\}} U(t,x)\ dx \le 0
\;\;\mbox{ for }\;\; t\ge t_0\,,
$$
whence
$$
\int_{\{ \vert x\vert\ge y(t)\}} U(t,x)\ dx \le \int_{\{ \vert
x\vert\ge \sigma(t_0)\}} U(t_0,x)\ dx = 0\,, \quad t\ge t_0\,.
$$
We have thus established that $\sigma(t)\le y(t)$ for $t\ge t_0$ from
which we readily conclude that 
\beqn
\label{f6}
\sigma(t) \le \sigma(t_0) + C\ \Vert U(t_0)\Vert_1^{q\xi(p-2)/p}\
(t-t_0)^{(1-N\xi(p-2))/p}\,, \quad t\ge t_0\,.
\eeqn

\smallskip

We next integrate (\ref{vhj1}) over $\RR^N$ and obtain
$$
\frac{d}{dt} \Vert U(t)\Vert_1 + \Vert\nabla U(t)\Vert_q^q = 0\,.
$$
Since the support of $U(t)$ is included in $B(0,\sigma(t))$, we infer
from Lemma~\ref{lef1} that
$$
\Vert\nabla U(t)\Vert_q^q = \int_{\{\vert x\vert<\sigma(t)} \vert\nabla
U(t,x)\vert^q\ dx \ge \frac{1}{\kappa\ \sigma(t)^{1/\xi}}\
\left( \int_{\{\vert x\vert < \sigma(t)\}} U(t,x)\ dx \right)^q =
\frac{ \Vert U(t)\Vert_1^q}{\kappa\ \sigma(t)^{1/\xi}}\,.
$$ 
Inserting this lower bound in the previous diferential equality gives
\beqn
\label{f7}
\frac{d}{dt} \Vert U(t)\Vert_1 +  \frac{1}{\kappa}\ \frac{\Vert
U(t)\Vert_1^q}{\sigma(t)^{1/\xi}} \le 0\,. 
\eeqn

Before going on we introduce the following notations: 
\bean
\Sigma(T) := \sup_{t\in [1,T]}{\left\{ t^{-A}\ \sigma(t) \right\}}\,,&
\quad & A:= \frac{q-p+1}{2q-p} \,, \\
L(T) := \sup_{t\in [1,T]}{\left\{ t^B\ \Vert U(t)\Vert_1 \right\}} \,,
& \quad & B:= \frac{(N+1)(q_*-q)}{2q-p}\,,  
\eean
for $T\ge 1$ and notice that $\Sigma(T)$ and $L(T)$ are well-defined
for each $T\ge 1$ while $A$ and $B$ satisfy
\beqn
\label{f8}
A + \frac{q\xi(p-2)}{p}\ B = \frac{1-N\xi(p-2)}{p} \;\;\mbox{ and
}\;\; 1 - \frac{A}{\xi} = (q-1)\ B \,.
\eeqn 
Fix $T\ge 1$. We infer from (\ref{f7}) that, if $t\in [1,T]$, 
\bean
& & \frac{d}{dt} \Vert U(t)\Vert_1 +  \frac{t^{-A/\xi}}{\kappa}\
\frac{\Vert U(t)\Vert_1^q}{t^{-A/\xi}\ \sigma(t)^{1/\xi}} \le 0\\
& & \frac{d}{dt} \Vert U(t)\Vert_1 +  \frac{1}{\kappa\
\Sigma(T)^{1/\xi}}\ \frac{\Vert U(t)\Vert_1^q}{t^{A/\xi}} \le 0\,,
\eean
which gives
\beqn
\label{f9}
\Vert U(t)\Vert_1 \le C\ \Sigma(T)^{1/((q-1)\xi)}\ \left( t^{(q-1)B} -
1 \right)^{-1/(q-1)}\,, \quad
t\in [1,T]\,,
\eeqn
after integration. Consider next $t\in [1,T]$. Either $t\le 4$ and it follows
from (\ref{f6}) with $t_0=1$ that 
$$
t^{-A}\ \sigma(t) \le t^{-A}\ \sigma(1) + C\ \Vert U(1)\Vert_1^{q\xi(p-2)/p}\
(t-1)^{(1-N\xi(p-2))/p}\ t^{-A} \le C\,.
$$
Or $t\ge 4$ and we infer from (\ref{f6}) with $t_0=t/2\ge 2$,
(\ref{f8}) and (\ref{f9}) 
that 
\bean
t^{-A}\ \sigma(t) & \le & t^{-A}\ \sigma\left( \frac{t}{2} \right) +
C\ \left\Vert U\left( \frac{t}{2} \right)\right\Vert_1^{q\xi(p-2)/p}\
t^{q\xi(p-2)B/p} \\
& \le & 2^{-A}\ \Sigma(T) + C\ \Sigma(T)^{(q(p-2))/(p(q-1))}\,.
\eean
Consequently, 
$$
t^{-A}\ \sigma(t) \le 2^{-A}\ \Sigma(T) + C\ \left( 1 +
\Sigma(T)^{(q(p-2))/(p(q-1))} \right)\,, \quad t\in [1,T]\,,
$$ 
from which we conclude that 
$$
\Sigma(T) \le 2^{-A}\ \Sigma(T) + C\ \left( 1 +
\Sigma(T)^{(q(p-2))/(p(q-1))} \right)\,.  
$$
Since $A>0$ and $q(p-2)<p(q-1)$ the above inequality entails that
$\Sigma(T)\le C$ for each $T\ge 1$, the constant $C$ being independent
of $T$. Recalling (\ref{f3}) we have thus proved that $\varrho(t)\le
\sigma(t)\le C\ t^A$ for $t\ge 1$, hence (\ref{rad3}). 

Furthermore the boundedness of $\Sigma(T)$ and (\ref{f9}) ensure that
$\Vert U(t)\Vert_1 \le C\ (t-1)^{-B}$ for $t\ge 1$ which, together
with (\ref{f2}), implies that 
\beqn
\label{f10}
\Vert u(t)\Vert_1\le C\ t^{-B}\,, \quad t\ge 2\,.
\eeqn
We have thus also established the assertion~(iii) of
Corollary~\ref{cora6}. \qed

\medskip

\noindent\textbf{Proof of Corollary~\ref{cora6}.} Assume first that
$q\in (1,p-1)$. Then, on the one hand, it follows from (\ref{rad1})
that there is $\varrho_\infty>0$ such that
$\varrho(t)\le\varrho_\infty$ for $t\ge 1$. On the other hand, we may
proceed as in the proof of (\ref{f7}) to establish that
\beqn
\label{f11}
\frac{d}{dt} \Vert u(t)\Vert_1 +  \frac{1}{\kappa}\ \frac{\Vert
u(t)\Vert_1^q}{\varrho(t)^{1/\xi}} \le 0\,.
\eeqn
Therefore, 
$$
\frac{d}{dt} \Vert u(t)\Vert_1 +  \frac{1}{\kappa}\ \frac{\Vert
u(t)\Vert_1^q}{\varrho_\infty^{1/\xi}} \le 0\,, \quad t\ge 1\,,
$$
from which (\ref{decL11}) readily follows. 

\smallskip

Similarly, if $q=p-1$, we infer from (\ref{rad2}) and (\ref{f11})
that, for $t\ge 2$,  
\bean
\Vert u(t)\Vert_1 & \le & C\ \left( \int_1^t (1+\ln{s})^{-1/\xi}\ ds
\right)^{-1/(q-1)} \\
& \le & C\ \left( \int_0^{\ln{t}} (1+s)^{-1/\xi}\ e^s\ ds
\right)^{-1/(q-1)} \\ 
& \le & C\ \left( (1+\ln{t})^{-1/\xi}\ (t-1) \right)^{-1/(q-1)}\,, 
\eean
which gives (\ref{decL12}). 

\smallskip

Since the case $q\in (p-1,q_*)$ has already been handled in the proof
of Theorem~\ref{tha5} (recall (\ref{f10})) we are left with the case
$q=q_*$. In that particular case, $\xi=\eta$ and we infer from
(\ref{rad4}) and (\ref{f11}) that 
$$
\frac{d}{dt} \Vert u(t)\Vert_1 +  \frac{C}{t}\ \Vert
u(t)\Vert_1^q \le 0\,, \quad t\ge 1\,,
$$
which gives (\ref{decL14}) by integration. \qed 


\section{Persistence of dead cores}\label{pode}

\noindent\textbf{Proof of Proposition~\ref{pra7}.} We first study
the one-dimensional case $N=1$. We consider a non-negative function
$y\in\mathcal{C}^1([0,\infty))$ to be specified later and proceed as
in the proof of Theorem~\ref{tha5} to deduce from (\ref{vhj1}) that 
\beqn
\label{pdc2}
\frac{d}{dt} \int_{-y(t)}^{y(t)} u(t,x)\ dx = \left[ \left(
\frac{p-1}{p-2}\ \left\vert \partial_x \left( u^{(p-2)/(p-1)}
\right)(t,x) \right\vert^{p-1} + y'(t) \right)\
u(t,x)\right]_{x=-y(t)}^{x=y(t)} 
\eeqn
On the one hand we infer from (\ref{ge1}) that
\bean
\frac{p-1}{p-2}\ \left\vert \partial_x \left( u^{(p-2)/(p-1)}
\right)(t,x) \right\vert^{p-1} & \le & \frac{p-1}{p-2}\ C(p,1)^{p-1}\
\Vert u_0\Vert_\infty^{(p-2)/p}\ t^{-(p-1)/p} \\ 
& \le & c_1\ \Vert u_0\Vert_\infty^{(p-2)/p}\ t^{-(p-1)/p}\,.
\eean
On the other hand, since $p-1>q$, we have
$\beta_{p,q}=\alpha_p=(p-2)/(p-1)$ and it follows from (\ref{ge2})
that  
\bean
\frac{p-1}{p-2}\ \left\vert \partial_x \left( u^{(p-2)/(p-1)}
\right)(t,x) \right\vert^{p-1} & \le & \frac{p-1}{p-2}\
C(p,q,1)^{p-1}\ \Vert u_0\Vert_\infty^{(p-1-q)/q}\ t^{-(p-1)/q} \\ 
& \le & c_2\ \Vert u_0\Vert_\infty^{(p-1-q)/q}\ t^{-(p-1)/q}\,.
\eean
Consequently, choosing
\beqn
\label{pdc4}
\left\{ 
\begin{array}{l}
y'(t) = - \min{\left\{ c_1\ \Vert u_0\Vert_\infty^{p-2}\ t^{-(p-1)/p}
, c_2\ \Vert u_0\Vert_\infty^{(p-1-q)/q}\ t^{-(p-1)/q} \right\}}\,, \\
\\ 
y(0)=R_0\,,
\end{array}
\right.
\eeqn
we have 
\beqn
\label{pdc3}
\frac{p-1}{p-2}\ \left\vert \partial_x \left( u^{(p-2)/(p-1)}
\right)(t,x) \right\vert^{p-1} \le - y'(t)\,. 
\eeqn
We then deduce from (\ref{pdc2}) and (\ref{pdc3}) that 
$$
\frac{d}{dt} \int_{-y(t)}^{y(t)} u(t,x)\ dx \le 0\,,
$$
whence
$$
\int_{-y(t)}^{y(t)} u(t,x)\ dx \le \int_{-R_0}^{R_0} u_0(x)\ dx = 0
\;\;\mbox{ for }\;\; t\ge 0\,. 
$$
Now it is actually possible to compute the function $y$ defined by
(\ref{pdc4}) and to see that  
$$
y(t) \ge y_\infty:=\lim_{s\to\infty} y(s) = R_0 - \delta_0\ \Vert
u_0\Vert_\infty^{(p-1-q)/(p-q)} 
$$
for some $\delta_0$ depending only on $c_1$, $c_2$, $p$, and $q$. Then
$u(t,x)=0$ for $x\in [-y_\infty,y_\infty]$ and $t\ge 0$, and
$y_\infty>0$ under the assumptions of Proposition~\ref{pra7}. 

\medskip

In several space dimensions $N\ge 2$, consider $\varepsilon\in
(0,R_0/2)$ and put 
$$
u_0^\varepsilon(x_1) := \left\{
\begin{array}{lcl}
\Vert u_0\Vert_\infty & \mbox{ if } & \vert x_1\vert\ge R_0\,, \\ \\
\displaystyle{\frac{\Vert u_0\Vert_\infty}{\varepsilon}\ (\vert
x_1\vert - R_0 + \varepsilon) } & \mbox{ if } & R_0 - \varepsilon \le
\vert x_1\vert\le R_0\,, \\ \\ 
0 & \mbox{ if } & \vert x_1\vert\le R_0 - \varepsilon\,, 
\end{array}
\right.
$$
Clearly, $u_0\le u_0^\varepsilon$ in $\RR^N$ and the comparison
principle entails that $u(t,x_1,x_2,\ldots,x_N)\le
u^\varepsilon(t,x_1)$ for $(t,x)\in [0,\infty)\times\RR^N$, where
$u^\varepsilon$ denotes the solution to (\ref{vhj1}) with initial
condition $u_0^\varepsilon$ and $N=1$. Choosing $\varepsilon$
appropriately small provides the expected result in the
$x_1$-direction. We proceed analogously in every direction to complete
the proof of Proposition~\ref{pra7}. \qed 


\appendix

\section{Proof of Lemma~\ref{leb1}}

\setcounter{thm}{0}
\setcounter{equation}{0}

Since $\partial_t u = \varphi'(v)\ \partial_t v$ and $\nabla u =
\varphi'(v)\ \nabla v$ we deduce from (\ref{b1}) that
$$
\partial_t v - a\ \Delta v - a\ \frac{\varphi''}{\varphi'}\ w
- 2\ a'\ \varphi'\ \varphi'' w^2 - 2\ a'\ {\varphi'}^2\
(\nabla v)^t D^2 v \nabla v + \frac{b'}{\varphi'} = 0. 
$$
Observing that
$$
(\nabla v)^t D^2 v \nabla v = \frac{1}{2}\ \nabla v\cdot \nabla w
\;\;\mbox{ and }\;\; \Delta w = 2\ \nabla v\cdot \nabla\Delta v + 2\
\sum_{i,j} \left\vert \partial_i \partial_j v \right\vert^2\,,
$$
elementary, but laborious calculation shows that
$$
\partial_t w - \mathcal{A}w + 2\ a\ \sum_{i,j} \left\vert \partial_i
\partial_j v \right\vert^2 + 2\ a'\ \varphi'\ \varphi''\ w\ \nabla
v\cdot \nabla w - \mathcal{V}\cdot\nabla w + 2\ \mathcal{S}_1\ w^2 +
2\ \mathcal{R}_2\ w = 0
$$
with
\beqn
\label{ap2}
\mathcal{S}_1 := - a\ \left( \frac{\varphi''}{\varphi'}
\right)' - 2\ a'\ \varphi'\ \varphi''\ \Delta v - 4\ a'' \left(
\varphi' \varphi'' \right)^2\ w^2 - 2\ a'\ w \left( 2 {\varphi''}^2 +
\varphi' \varphi''' \right)\,,
\eeqn
and 
\bear
\nonumber
\mathcal{V} & := & 2\ \left[ a\ \frac{\varphi''}{\varphi'} + a'\
{\varphi'}^2\ \left( \Delta v + \frac{2 \varphi''}{\varphi'}\ w
\right) \right]\ \nabla v \\
\nonumber
& + & 4\ \varphi'\ \varphi'' \left[ \left( a''\ {\varphi'}^2\ w + 3\ a'
\right) + a''\ {\varphi'}^2\ w \right]\ w\ \nabla v \\
\label{ap1}
& + & 2\ \left[ a''\ {\varphi'}^4\ \nabla v\cdot \nabla w - b'\
\varphi' \right]\ \nabla v +a'\ {\varphi'}^2\ \nabla w\,.
\eear

In order to handle the term involving $\Delta v$ in $\mathcal{S}_1$ we
proceed as in \cite{Be81}: more precisely we have
\bean
& & 2\ a\ \sum_{i,j} \left\vert \partial_i \partial_j v \right\vert^2 + 2\ a'\
\varphi'\ \varphi''\ w\ \nabla v\cdot \nabla w - 4\ a'\ \varphi'\
\varphi''\ \Delta v\ w^2 \\
& = & 4\ a'\ \varphi'\ \varphi''\ w \left( \frac{1}{2}\ \nabla
v\cdot\nabla w - w\ \Delta v \right) +  2\ a\ \sum_{i,j} \left\vert \partial_i
\partial_j v \right\vert^2 \\
& = & 4\ a'\ \varphi'\ \varphi''\ w \left( \sum_{i,j} \partial_i
\partial_j v\ \partial_i v\ \partial_j v - w\ \sum_i \partial_i^2 v
\right) +  2\ a\ \sum_{i,j} \left\vert \partial_i \partial_j v
\right\vert^2 \\
& = & \sum_i \left\{ 2\ a\ \left\vert\partial_i^2 v\right\vert^2 + 4\
a'\ \varphi'\ \varphi''\ w\ \left(\left\vert\partial_i v\right\vert^2 -
w \right)\ \partial_i^2 v \right\} \\
& + & \sum_{i\ne j} \left\{ 2\ a\ \left\vert\partial_i \partial_j
v\right\vert^2 + 4\ 
a'\ \varphi'\ \varphi''\ w \partial_i \partial_j v\ \partial_i v\
\partial_j v \right\} \\ 
& = & 2\ a\ \sum_{i} \left\{ \partial_i^2 v + \frac{a'}{a}\ \varphi'\
\varphi''\ w \left(\left\vert\partial_i v\right\vert^2 - w \right)
\right\}^2 \\
& - & 2\ \sum_{i} \frac{{a'}^2}{a}\ \left( \varphi' \varphi''
\right)^2\ w^2\ \left(\left\vert\partial_i v\right\vert^2 - w
\right)^2 \\
& + & 2\ a\ \sum_{i\ne j} \left\{ \partial_i \partial_j v +
\frac{a'}{a}\ \varphi'\ \varphi''\ w\ \partial_i v\ \partial_j v \right\}^2 \\
& - & 2\ \sum_{i\ne j} \frac{{a'}^2}{a}\ \left( \varphi' \varphi''
\right)^2\ w^2\ \left\vert\partial_i v\right\vert^2\
\left\vert\partial_j v\right\vert^2 \\
& \ge & -2\ (N-1)\ \frac{{a'}^2}{a}\ \left( \varphi' \varphi''
\right)^2\ w^2\,.
\eean
Consequently,
$$
2\ a\ \sum_{i,j} \left\vert \partial_i \partial_j v
\right\vert^2 + 2\ a'\ \varphi'\ \varphi''\ w\ \nabla v\cdot \nabla w
+ 2\ \mathcal{S}_1\ w^2\ge 2\ \mathcal{R}_1\ w^2\,,
$$
which completes the proof of the first assertion of Lemma~\ref{leb1}. 

\medskip

In the case where $x\longmapsto u(t,x)$ is radially symmetric and
non-increasing for each $t\ge 0$, we have $u(t,x)=U(t,|x|)$ for
$(t,x)\in [0,\infty)\times\RR^N$ and $\partial_r U(t,r)\le 0$ for
$(t,r)\in [0,\infty)\times [0,\infty)$. Introducing
$V=\varphi^{-1}(U)$ we have $v(t,x)=V(t,|x|)$ and the monotonicity of
$\varphi$ warrants that $\partial_r V(t,r)\le 0$. In addition, owing
to the non-negativity of $a'$, $\varphi'$ and $\varphi''$, we have
\bean
& & 2\ a'\ \varphi'\ \varphi''\ w\ \nabla v\cdot \nabla w - 4\ a'\
\varphi'\ \varphi''\ w^2\ \Delta v\\
& = & 2\ a'\ \varphi'\ \varphi''\ w\ \left[ 2\ \left\vert \partial_r
V\right\vert^2\ \partial_r^2 V  - 2\ \left\vert \partial_r
V\right\vert^2\ \left( \partial_r^2 V + \frac{N-1}{r}\ \partial_r V
\right) \right] \\ 
& \ge & 0\,,
\eean 
from which we deduce that
$$ 
2\ a'\ \varphi'\ \varphi''\ w\ \nabla v\cdot \nabla w + 2\
\mathcal{S}_1\ w^2 \ge 2\ \mathcal{R}_1^r\ w^2\,,
$$
and end the proof of Lemma~\ref{leb1}.  \qed


\end{document}